\documentclass[12pt,psamsfonts]{article}
\usepackage{amsmath}
\usepackage{amssymb}
\newtheorem{theo}{Theorem}[section]
\newtheorem{remarkk}[theo]{Remark}
\newenvironment{rem}{\begin{remarkk}\rm}{\end{remarkk}}
\newtheorem{definition}[theo]{Definition}
\newenvironment{devery easyfi}{\begin{definition}\rm}{\end{definition}}
\newtheorem{prop}[theo] {Proposition}
\newtheorem{cor}[theo]{Corollary}
\newtheorem{lemma}[theo]{Lemma}
\newtheorem{example}[theo]{Example}

\newcommand{\CC}{\ensuremath{\mathbb{C}}}

\newcommand{\ZZ}{\ensuremath{\mathbb{Z}}}

\newcommand{\M}{\ensuremath{\mathbb{M}}}
\newcommand{\LL}{\ensuremath{\mathbb{L}}}

\newcommand{\hol}{\ensuremath{\mathcal{O}}}

\newcommand{\PP}{\ensuremath{\mathbb{P}}}

\newcommand{\X}{\ensuremath{\mathcal{X}}}

\newcommand{\ra}{\ensuremath{\rightarrow}}

\newcommand{\G}{\ensuremath{\mathcal{G}}}
\newcommand{\F}{\ensuremath{\mathcal{F}}}

\def\Bbb{\bf}
\def\P{{\Bbb P}}

\def\Z{{\Bbb Z}}

\newcommand{\qed}{\hspace*{\fill}$Q.E.D.$}
\begin{document}
\title{Symmetry and Variation of Hodge Structures}
\author{I. C. Bauer, F. Catanese \\
\small Mathematisches Institut der Universit\"at Bayreuth \\
\small Universit\"atsstr. 30,
\small D-95447 Bayreuth\\
    This article is dedicated to Yum Tong Siu's
   60-th birthday.
}
\date{\today}

\maketitle

\section{Introduction}

The Torelli theorem (cf.\cite{tor}, \cite{andr}, \cite{weil1}) states
that two algebraic curves are
isomorphic if and only if their  Jacobian varieties  are isomorphic
as polarized Abelian varieties.

Andre' Weil (\cite{weil2}) set up a program for doing arithmetics on
K3 surfaces, based on a Torelli
type theorem, which was later proven through the effort of several authors (for this long
story and related references we refer to Chapter VIII of the book
\cite{BPV}). This result is crucial
for answering questions about the existence of K3 surfaces or
families thereof possessing certain curve
configurations.

The general Torelli question, as set up by Griffiths (\cite{grif68}, \cite{grif70}, cf. also
\cite{gr-s}, \cite{grif84}), is to associate to each projective
variety $X$ of general type its Hodge
structure of weight $n = dim (X)$, and ask whether the corresponding
"period map"
$\psi_n$ is injective on the local moduli space (or Kuranishi space of $X$).

It was known since long time that, as soon as the dimension is at
least two, there are families of
varieties without Hodge structures, which are not rigid.
Surfaces of general type with $q=p_g=0$
were constructed in the 30's by Campedelli and Godeaux (\cite{Cam},
\cite{god}), and for instance, in
the case of the Godeaux surfaces, the Kuranishi family has dimension $8$.

A natural question which arises is: under which hypothesis on $X$ is
a local Torelli theorem valid for
the Hodge structure of weight $n=dim X$?

 In other words, when is the
local period map $\psi_n$ a local
embedding?

The question is already quite open in dimension $n=2$, and the
hypothesis of requiring $X$ to be simply
connected only made the search of counterexamples more complicated (cf. e.g. \cite{cat} for a brief
account of these examples).

On the other hand, Lieberman, Peters and Wilsker (cf. \cite{l-p-w}) made
clear that, thanks to Griffiths'
interpretation of the derivative of the period map as a
cup-product (cf. further below in the
introduction), the infinitesimal Torelli theorem (injectivity of the
derivative) would follow from the
vanishing of a certain Koszul cohomology group.

This approach was later developed by several authors (cf. e.g.
\cite{green1}, \cite{green2}, \cite{flenner}, \cite{cox}) who essentially proved that, given a
construction involving some degree
$d$ (taking hypersurfaces of degree $d$ in some fixed manifold, or
complete intersections), then
infinitesimal Torelli holds for $d$ sufficiently large. In some
sense,  results  of this kind parallel
Serre's vanishing theorem $B(n)$, and it would be beautiful to give
precise geometrical conditions
which would ensure the validity of the infinitesimal Torelli theorem
e.g. for varieties of general type
with ample canonical bundle.

Certainly it has been up to now an open question whether the
condition that the canonical system be
very ample is sufficient for this purpose (the condition was not
holding for all hitherto known
counterexamples). Note that, if the canonical bundle is ample but not
very ample, the infinitesimal
Torelli theorem already fails in dimension $1$, as the case of
hyperelliptic curves already shows.

Unfortunately, we show in this paper that this condition might not be
sufficient, indeed we provide a series of
examples of surfaces of general type for which  the infinitesimal
Torelli theorem fails on the whole
moduli space in the worst possible way, namely, the period map has
all fibres of positive
dimension. This holds in spite of the fact that the canonical system
is generally quasi very ample (i.e., we prove that it is a birational 
 morphism which yields a local embedding on the complement 
of a finite set).

Our examples do not rule out the possibility that for very ample canonical
system the period map may be generically finite. Indeed, in our examples
one recovers the "missing" Hodge structure from the geometry of the
singular locus of the canonical image.

We raise therefore the question: what  geometric properties
are required, for a variety with
very ample canonical system, for the validity of the infinitesimal
Torelli theorem? For instance, in
terms of the geometry of the canonical image?

    We proceed now to a more detailed information, introducing the
standard notation and our situation.

Let $X$ be a smooth algebraic variety over the complex numbers, which
for simplicity we assume to have
ample canonical bundle $K_X$. Then $K_X$ defines a natural
polarization on $X$ and we know from the
Kodaira - Spencer - Kuranishi theory (cf. \cite{K-M}, \cite{K}) that
there exists a semiuniversal
deformation $p : \X \longrightarrow (Y,y_0)$ of $X$. In particular,
the tangent space $T_{Y,y_0}$ of $Y$
in $y_0$ is naturally isomorphic to $H^1(X,T_X)$ and the dimension of
$Y$ in $y_0$ is at least
$dim  \ H^1(X,T_X) - dim \ H^2(X,T_X)$. \\

For each $k \in \{1,\ldots,dimX\}$ we have a corresponding variation
of Hodge structure
$$ (H_{\mathbb{Z}} =
\mathcal{R}^kp_*(\mathbb{Z}),H^{p,q}(y),Q),
$$
where $p+q = k$, $H_{\mathbb{Z}} \otimes \mathbb{C} =
\bigoplus_{p+q=k}H^{p,q}(y)$, and the polarization $Q$ is a
quadratic form on $H_{\mathbb{Z}}$
for which the subspaces $H^{p,q}(y)$ are pairwise orthogonal. To this
variation of Hodge structure
there is associated a holomorphic map $\Phi : Y \longrightarrow D$,
where $D$ is the classifying domain
of polarized Hodge structures of type $(h^{k,0}, h^{k-1,1}, \ldots
,h^{0,k})$  and $h^{p,q} =
dimH^{p,q}(y)$ (cf. e.g. \cite{gr-s}).\\

The {\em infinitesimal Torelli theorem} is said to hold for $X$ if
and only if the differential $d\Phi$
of
$\Phi$ is injective on $T_{Y,y_0} $.\\

In the late sixties P. Griffiths (\cite{grif68}, \cite{grif70}) showed
that the differential $d\Phi$ of
the period map is given by the map below, induced by cup product

$$ d\Phi : H^1(X,T_X) \longrightarrow \bigoplus_{p+q=k}
Hom(H^{p,q}(X), H^{p-1,q+1}(X)).
$$

Classically, for a smooth curve $C$ the infinitesimal Torelli theorem
holds iff $g(C) = 1,2$ or iff
$g(C) \geq 3$ and $C$ is not hyperelliptic. \\

P. Griffiths posed the problem to determine the class of varieties with
$K_X$ (sufficiently) ample for which infinitesimal Torelli holds.
Even if there are several
counterexamples to the infinitesimal Torelli theorem known, there is
still the hope that for projective
manifolds with sufficiently ample canonical bundle the result should
be true. \\

The simple underlying idea of this paper is the following: assume
that $X$ is the quotient of a
smooth variety $Z$ by the free action of a finite group $G$.
   Then, in the Griffiths cup product map, if we replace $Z$ by $X$, we
must replace each cohomology group
by the subspace of $G$-invariants.

However, even if infinitesimally Torelli holds for $Z$, it does not
need to hold for $X$ any longer, by
the simple algebraic observation that the tensor product of
invariants is much smaller than
the subspace of invariants in the tensor product of $G$-representations.

A similar philosophy was used by S. Usui (\cite{Usui}) to justify
failure of infinitesimal Torelli
for varieties $Z$ with automorphisms, and in \cite{cat89} to produce
generalized examples of
everywhere non reduced moduli spaces (in the latter paper one needed
the action not to be
free, whereas in the former one had no restriction, but the "bad"
varieties $Z$ with special
automorphisms were not generic in the moduli space).

So, we shall study here the situation for quotients of projective
manifolds  by a  free action of a
finite group and it will turn out that there are examples of surfaces
having quasi very ample canonical
bundle and for which infinitesimal Torelli fails. \\

The most natural candidates, however, do not suffice for our
purposes: quotients of a hypersurface
   $ Z \subset \PP^3$ of degree $d$ by the action of a finite group $G$
of order $m$.
The first example is given by the classical Godeaux surfaces
($d=m=5$), for which however there is no
Hodge structure, and the second simplest example, with $d=m=7$,
satisfies already local Torelli
at the general point (it does so at the quotient of the Fermat surface).

For this reason we have to resort to quotients of products of curves,
the so called surfaces isogenous
to a (higher) product, whose moduli spaces were thoroughly 
investigated in \cite{cat00}
  (cf. also \cite{cat03}).

The simplest examples will be the ones where the group $G$ is cyclic.
When $G$ has order two, we will
get positive dimensional fibres of the period map which are relatives
of the positive dimensional
fibres of the Prym map (cf. \cite{naranjo}). This case however is not
fully satisfactory because we
want $K_X$ to be quasi very ample.

\bigskip

The following
will be one of our main results.

\begin{theo} For any natural number $k \geq 2$ there exists a family
of surfaces $\mathcal{S}_k$ of
dimension $3k+2$ such that the following hold:\\ 1) for each
$S_k \in \mathcal{S}_k$ the
infinitesimal Torelli map
$$ d\Phi_2 : H^1(S_k,T_{S_k}) \longrightarrow Hom(H^{2,0}(S_k), H^{1,1}(S_k))
$$ has a kernel of dimension at least two, i.e., the period map
$\Phi_2$ belonging to the Hodge
structures of weight $2$ has fibres  of dimension at least two  over
each point. \\ 2) for a general
$S_k \in \mathcal{S}_k$ the canonical divisor $K_{S_k}$ is
quasi very ample, i.e., it gives a birational  morphism which yields a
local embedding on the complement of a finite set;
\\ 3)
$h^1(S_k,T_{S_k}) = 3k+2$; in particular $\mathcal{S}_k$ is a 
generically smooth
irreducible component in the moduli
space.
\end{theo}

We emphasise once more here the important fact that the surfaces we 
are considering all
have  unobstructed local moduli spaces (i.e.,
the basis of the Kuranishi family is smooth, and indeed an open set
in $ H^1(S,T_{S})$).

After preliminaries in section $1$, the above main theorem will be
shown in sections $2$ and $3$.

In section $4$, however, we just take the opposite point of view and
ask the converse question:
since the Torelli theorem holds for curves, may it be that one can
completely characterize the case
where a global $n$-tuple Torelli theorem does not hold for varieties
isogenous to a product of curves?

Here,  $n$-tuple Torelli means that we can reconstruct the variety
from all of its Hodge structure
(i.e., we do not restrict only to weight $n = dim (X)$).

The question looks interesting in general, but rather complicated
already for surfaces.
For this reason we content ourself with giving sufficient conditions
for the validity of double Torelli
for surfaces isogenous to a product (theorem 4.4). In proving this
result, we establish some
intermediate technical results  which may be of
independent interest (lemma 4.6.,
proposition 4.5). These results concern the eigenspaces, for the 
action of a cyclic
group on an algebraic curve $C$, inside the space of
holomorphic differentials on $C$.

Finally, in section $5$, we discuss another series of examples
where we are not able to show global
double Torelli, but where the local period map is injective.

\section{Symmetry imposes failure of Torelli type theorems}

As already mentioned, in this paper we  exploit the above main idea
as a simple tool for providing
generic counterexamples to Torelli theorems. A similar observation was
made (in a different direction,
however, cf.
\cite{Usui}) to show that often the failure of Torelli theorems (especially
infinitesimal Torelli) is due to the presence of symmetry on the
variety under consideration.\\

Our situation is more specific: we assume that $G$ is a finite  group
acting freely on a smooth
algebraic variety $Z$ of dimension $n$. Let $X$ be the quotient $Z / G$.

Then the infinitesimal  Torelli theorem for the priods of $n$-forms holds for 
$Z$ if and only if we have surjectivity of

$$ H^{n-1}(Z, \Omega^1_Z) \otimes H^0(Z, \Omega^n_Z) \longrightarrow
H^{n-1} (Z, \Omega ^1_Z \otimes
\Omega^{n-1}_Z).
$$

Whence, if infinitesimal Torelli holds for $Z$, we have surjectivity
also of the map
$$ 
(H^{n-1}(Z, \Omega^1_Z) \otimes H^0(Z, \Omega^n_Z))^G  \rightarrow
(H^{n-1} (Z, \Omega ^1_Z \otimes
\Omega^{n-1}_Z))^G = H^{n-1} (X, \Omega ^1_X \otimes \Omega^{n-1}_X).
$$

   For simplicity of notation assume that $G$ is abelian (a quite
similar fact holds in the general case,
we shall however concentrate ourselves on the more tractable case
where $G$ is abelian):  by Schur's
lemma we have
$$ 
(H^{n-1}(Z, \Omega^1_Z) \otimes H^0(Z, \Omega^n_Z))^G =
\bigoplus_{\chi \in G^*} H^{n-1}(Z,
\Omega^1_Z)^{\chi} \otimes H^0(Z, \Omega^n_Z)^{\chi^*}.
$$

In order to get counterexamples for the infinitesimal 
Torelli theorem for periods of $n$-forms, we have to look for
situations where, although
$$
\bigoplus_{\chi \in G^*} H^{n-1}(Z, \Omega^1_Z)^{\chi} \otimes H^0(Z,
\Omega^n_Z)^{\chi^*} \rightarrow
(H^{n-1} (Z, \Omega ^1_Z \otimes \Omega^{n-1}_Z))^G
$$
is surjective, still
$$ H^{n-1}(Z, \Omega^1_Z)^G \otimes H^0(Z, \Omega^n_Z)^G \rightarrow
H^{n-1} (Z, \Omega ^1_Z \otimes
\Omega^{n-1}_Z)^G
$$
fails to be surjective. \\

The first case where this occurs is the one of a classical Godeaux surface.
This is the quotient of a  smooth
quintic $Z$ in $\mathbb{P}^3$ on which $G := \mathbb{Z}/5$ acts freely.
In this case however $H
^0(\Omega^2_Z)^G = 0$, whence there is no period map. \\

We tried to cook up other examples as quotients, e.g. of complete
intersections in projective space, but
the examples which worked best were the cases of quotients of
products of curves. For these we will
produce examples where the canonical bundle  remains ample or even
quasi very ample, but infinitesimal Torelli
fails.

\section{Surfaces isogenous to a product}

Let us recall the notion of surfaces isogenous to a higher product
(prop. 3.11 of \cite{cat00}
ensures that the following two properties 1) and 2) of a surface are
equivalent).

\begin{definition}  A surface $S$ is said to be {\em isogenous to a
higher  product} if and only if,
equivalently, either

1)  $S$ admits a finite unramified covering which is isomorphic to a
product of curves of genera at
least two, or \\ 2) $S$ is a quotient  $S := (C_1 \times C_2) / G$,
where the $C_i$'s are  curves of
genus least two, and $G$ is a finite group acting freely on $Z:= (C_1
\times C_2)$.

We have two cases: the {\em mixed case}, where the action of $G$
exchanges the two factors (and
then $C_1$, $C_2$ are isomorphic), and the {\em unmixed case}
where $G$ acts via a product action.
\end{definition}

We shall assume throughout that we have such a surface and that we
are in the unmixed case, thus we
have a finite group $G$ acting on two curves $C_1$, $C_2$  with genera
$g_1$, $g_2 \geq 2$, and  acting freely by the product action on
$Z := C_1 \times C_2$. We will now examine the infinitesimal Torelli
map of the quotient $S := C_1
\times C_2 /G.$

In the following theorem we shall use the standard notation: given a  finite
group
$G$, we let
$G^*$ be the set of characters of 
irreducible representations of
$G$, and, for $\chi \in G^*$ and $V$
a $G$-representation, we denote by $V^{\chi}$ the
$\chi$-isotypical component. Finally, for $\chi \in G^*$,  $\chi^*$
denotes the character of the dual
irreducible representation.

\begin{theo} \label{thm1} Let $G$, $C_1$, $C_2$ and $S$ be as above, and assume that
$G$ is abelian.

  Then
    the infinitesimal Torelli theorem holds for $S$, i.e., we have the
surjectivity of the following
linear map
$$ 
d\Phi_2^* : H^1(S,\Omega^1_S) \otimes H^0(S,\Omega^2_S)
 \longrightarrow H^1(S,\Omega^1_S \otimes \Omega^2_S),
$$
if and only if the two maps
$$
\bigoplus_{\chi \in G^* : H^0(C_2,\Omega^1_{C_2})^{\chi} \ \neq 0} H^0(C_1,
\Omega^1_{C_1})^{\chi} \otimes H^0(C_1, \Omega^1_{C_1})^{\chi^*}
\longrightarrow H^0(C_1,\Omega^1_{C_1} \otimes \Omega^1_{C_1})^G
$$
and
$$
\bigoplus_{\chi \in G^* : H^0(C_1,\Omega^1_{C_1})^{\chi} \ \neq 0} H^0(C_2,
\Omega^1_{C_2})^{\chi} \otimes H^0(C_2, \Omega^1_{C_2})^{\chi^*}
\longrightarrow H^0(C_2,\Omega^1_{C_2} \otimes \Omega^1_{C_2})^G
$$
are both surjective.
\end{theo}

{\em Proof.} \\ Let $Z := C_1 \times C_2$ and  denote the two
projections to $C_1$ resp. $C_2$ by
$p_1$, resp. $p_2$. Then $\Omega^1_Z = p_1^*
\Omega^1_{C_1} \oplus p_2^* \Omega^1_{C_2}$. Therefore we obtain

$$ H^1(S, \Omega^1_S) = H^1(Z,\Omega^1_Z)^G = H^1(Z, (\Omega^1_{C_1}
\boxtimes \mathcal{O}_{C_2})  \oplus  (\mathcal{O}_{C_1}  \boxtimes
\Omega^1_{C_2}))^G =
$$
$$
(H^1(C_1,\Omega^1_{C_1}) \otimes H^0(C_2, \mathcal{O}_{C_2}))
\oplus (H^0(C_1,\mathcal{O}_{C_1})
\otimes H^1(C_2, \Omega^1_{C_2})) \oplus
$$
$$
\bigoplus_{\chi \in G^*} (H^0(C_1,\Omega^1_{C_1})^{\chi}  \otimes H^1(C_2,
\mathcal{O}_{C_2})^{\chi^*}) \bigoplus_{\chi \in G^*} (
H^1(C_1,\mathcal{O}_{C_1})^{\chi} \otimes
H^0(C_2,
\Omega^1_{C_2})^{\chi^*}).
$$

Here, we use the convenient notation for the external tensor product
of coherent sheaves
$$
\F \boxtimes \G = p_1^*\F \otimes p_2^* \G.
$$
We remark that the above holds by  K\"unneth's formula and Schur's
lemma. Moreover, we use that
$H^1(C_1,\Omega^1_{C_1})$ and
$H^1(C_2,\Omega^1_{C_2})$ are automatically invariant, since every
automorphism of $C_i$ maps the
fundamental class to itself.

Again, using Schur's lemma we get
$$ H^0(S,\Omega^2_S) = H^0(Z,\Omega^2_Z)^G = H^0(Z,\Omega^1_{C_1} \boxtimes
\Omega^1_{C_2})^G =
$$
$$
\bigoplus_{\chi \in G^*} H^0(C_1, \Omega^1_{C_1})^{\chi} \otimes H^0(C_2,
\Omega^1_{C_2})^{\chi^*}.
$$
Moreover, we have
$$ H^1(S, \Omega^1_S \otimes \Omega^2_S) = H^1(Z, \Omega^1_Z \otimes
\Omega^2_Z)^G =
$$
$$ H^1(Z, (\Omega^1_{C_1})^2 \boxtimes \Omega^1_{C_2})^G \oplus H^1(Z,
\Omega^1_{C_1} \boxtimes (\Omega^1_{C_2})^2)^G =
$$
$$ (H^1(C_1, (\Omega^1_{C_1})^2) \otimes H^0(C_2, \Omega^1_{C_2}))^G
\oplus (H^0(C_1,
(\Omega^1_{C_1})^2) \otimes H^1(C_2, \Omega^1_{C_2}))^G
$$
$$
\oplus ( H^1(C_1, \Omega^1_{C_1}) \otimes H^0(C_2, (\Omega^1_{C_2})^2))^G
\oplus (H^0(C_1, \Omega^1_{C_1}) \otimes H^1(C_2, (\Omega^1_{C_2})^2))^G.
$$

Using now the fact that
$$ H^1(C_i, (\Omega^1_{C_i})^2) = 0,~~ H^1(C_i, \Omega^1_{C_i}) =
H^1(C_i, \Omega^1_{C_i})^G \ \ \ for \ \  i=1,2
$$
we obtain the simpler expression
$$ 
H^1(S, \Omega^1_S \otimes \Omega^2_S) =
$$
$$
(H^0(C_1, \mathcal{O}_{C_1}(2K_{C_1}))^G
\otimes H^1(C_2, \mathcal{O}_{C_2}(K_{C_2})) \oplus
$$
$$
\oplus (H^1(C_1, \mathcal{O}_{C_1}(K_{C_1}))
\otimes H^0(C_2, \mathcal{O}_{C_2}(2K_{C_2}))^G).
$$

By the non-degeneracy of the Serre duality,
$ \forall {\chi \in G^*}$ such that $H^0(C_i, \Omega^1_{C_i})^{\chi}  \neq 0$,
$H^0(C_i, \Omega^1_{C_i})^{\chi} \otimes H^1(C_i, \hol_{C_i})^{\chi^*}
\rightarrow H^1(C_i, \Omega^1_{C_i})$ is onto, whence we conclude that the dual Torelli map
$$ 
d\Phi^* : H^1(S,\Omega^1_S) \otimes H^0(S,\Omega^2_S)
\longrightarrow H^1(S,\Omega^1_S \otimes
\Omega^2_S)
$$
is surjective if and only if the two maps
$$
\bigoplus_{\chi \in G^* : H^0(C_2,\Omega^1_{C_2})^{\chi} \neq 0} H^0(C_1,
\Omega^1_{C_1})^{\chi} \otimes H^0(C_1, \Omega^1_{C_1})^{\chi^*}
\longrightarrow H^0(C_1,\Omega^1_{C_1} \otimes \Omega^1_{C_1})^G
$$
and
$$
\bigoplus_{\chi \in G^* : H^0(C_1,\Omega^1_{C_1})^{\chi} \neq 0} H^0(C_2,
\Omega^1_{C_2})^{\chi} \otimes H^0(C_2, \Omega^1_{C_2})^{\chi^*}
\longrightarrow H^0(C_2,\Omega^1_{C_2} \otimes \Omega^1_{C_2})^G
$$
are both surjective.

\qed

\bigskip We assume  now  that $G$ is a cyclic group of order $d$.
    Let $C_1$ and $C_2$ be smooth curves such that \\

1) $G$ acts freely on $C_1$,  \\

2) $G$ acts on $C_2$ in such a way that $C_2/G \cong  \mathbb{P}^1$. \\

Then $G$ acts freely on $Z := C_1 \times C_2$.

\begin{prop} \label{torellifalse} Assume that $C_1, C_2$ are as above, and
let $C_1':= C_1/G$.
\\  For $G = \mathbb{Z}/d$ the infinitesimal Torelli theorem does not hold for
$S := (C_1 \times C_2)/G$, if, $g'$ being the genus of $C_1'$,
\begin{itemize}
\item
$d = 2$ and $ 2 \leq g' \leq 5$
\item
$d=3$ and  $  g' = 3$
\item
$ 3 \leq d \leq 5$ and $g' = 2$.

\end{itemize}
\end{prop}

{\em Proof.} As we have seen in theorem \ref{thm1} it suffices to show that
$$
\Phi : \bigoplus_{\chi \in G^* : H^0(C_2,\Omega^1_{C_2})^{\chi} \ \neq 0}
    H^0(C_1, \Omega^1_{C_1})^{\chi} \otimes H^0(C_1, \Omega^1_{C_1})^{\chi^*}
\longrightarrow H^0(C_1,\Omega^1_{C_1} \otimes \Omega^1_{C_1})^G
$$
is not surjective. \\

Let us fix an isomorphism of $G$ with the group of $d$-th roots of
unity, and let us then denote by
$H^0(C_1,\mathcal{O}_{C_1}(K_{C_1}))^i$ the eigenspace of
$H^0(C_1,\mathcal{O}_{C_1}(K_{C_1}))$ belonging to the character $ i
\in (\Z / d)$. We have
$$ H^0(C_1,\mathcal{O}_{C_1}(K_{C_1}))^i =
H^0(C_1',\mathcal{O}_{C_1'}(K_{C_1'}) \otimes \mathcal{L}^i),
$$
where $\mathcal{L} \in Pic(C_1')$ is an element of precisely $d$ -
torsion, whence
$$ {\rm dim} \ H^0(C_1,\mathcal{O}_{C_1}(K_{C_1}))^i = g' - 1.$$
Therefore we have
$$ {\rm dim} \ \Phi \ (\bigoplus_{\chi \in G^* :
H^0(C_2,\Omega^1_{C_2})^{\chi} \neq 0} H^0(C_1,
\Omega^1_{C_1})^{\chi} \otimes H^0(C_1, \Omega^1_{C_1})^{\chi^*}) =
$$
$$ {\rm dim} \ \Phi \ (\bigoplus_{\chi \neq 0} H^0(C_1,
\Omega^1_{C_1})^{\chi} \otimes H^0(C_1, \Omega^1_{C_1})^{\chi^*})
\leq (\frac{d}{2} - 1)(g' - 1)^2 +
\frac{g'(g' - 1)}{2},
$$
for $d$ even and
$$ {\rm dim} \  \Phi \ (\bigoplus_{\chi \neq 0} H^0(C_1,
\Omega^1_{C_1})^{\chi} \otimes H^0(C_1, \Omega^1_{C_1})^{\chi^*})
\leq \frac{d-1}{2}(g' - 1)^2,
$$
for $d$ odd. \\ 
Observe now that
$$ {\rm dim } \ H^0(C_1, \mathcal{O}_{C_1}(2K_{C_1}))^G = h^0(C_1',
\mathcal{O}_{C_1'}(2K_{C_1'})) = 3g' - 3.
$$

From this we see immediately that
$$ {\rm dim} \ \Phi \ (\bigoplus_{\chi \neq 0} H^0(C_1,
\Omega^1_{C_1})^{\chi} \otimes H^0(C_1, \Omega^1_{C_1})^{\chi^*}) <
{\rm dim} \  H^0(C_1,
\mathcal{O}_{C_1}(2K_{C_1}))^G
$$
for
\begin{itemize}
\item
$d = 2$ and $ 2 \leq g' \leq 5$
\item
$d=3$ and  $  g' = 3$
\item
$ 3 \leq d \leq 5$ and $g' = 2$.
\end{itemize}
    \qed

\bigskip

We see here that the infinitesimal Torelli theorem fails simply for
reasons of dimension. Also it is
obvious from the formulae that we cannot get such an easy failure as
soon as $g'$ and $d$ become bigger.

\begin{rem} We have now constructed a series of counterexamples to
the infinitesimal Torelli theorem.
Of course we are now interested to see whether it is possible to
obtain that $K_S$ be very ample. If
$d=2$ there is no hope, since $\Z / 2$ induces on $C_2$ the
hyperelliptic involution, so that the
canonical map of $S$ is of degree $ \geq 2$: whence we  obtain that
$K_S$ is ample, but not very ample, or quasi very ample. Therefore we will
concentrate on the second and third case.
\end{rem}

Before constructing a concrete family of examples, we recall some of
the notations and results from
Pardini's article on abelian covers (\cite{Pardini1}): \\

Let $\pi : X \longrightarrow Y$ be a (finite) abelian cover with
group $G$, where $Y$ is smooth and $X$
is a normal variety. Then $\pi$ is flat and the action of $G$ induces
a splitting
$$
\pi_* \mathcal{O}_X = \bigoplus_{\chi \in G^*} L^{-1}_{\chi},
$$
where $G$ acts on $L^{-1}_{\chi}$ via the character $\chi$. The
invariant summand $L_0$ is isomorphic
to $\mathcal{O}_Y$. We denote by $D$ the branching divisor of $\pi$.
\\ We denote by $\mathfrak{S}$ the
set of cyclic subgroups of $G$ and for each $H \in \mathfrak{S}$ we denote by
$S_H$ the set of generators of the group of characters $H^*$. Then we can write
$$ D = \sum_{H \in \mathfrak{S}} \sum_{\psi \in S_H} D_{H, \psi},
$$
where $D_{H, \psi}$ is the sum of all components of $D$ having
inertia group $H$ and character $\psi$.
\\ We recall that the inertia group $H$ of a component $T$ of the
ramification divisor $R$ is defined by
$$ H = \{ h \in G : hx = x,~ \forall x \in T \}.
$$
Moreover, we associate to $T$ a generator $\psi_T$ of $H^*$, in fact
 there is a parameter $t$
of $\mathcal{O}_{X,T}$ such that the action of $H$ is given by
$$ ht = \psi_T(h)t, ~~ \forall h \in H.
$$
Now, if $S$ is a component of $D$, then all the components of
$\pi^{-1} (S)$ have the same inertia
group and isomorphic representations on the cotangent space at the
corresponding points of $X$ (since
$G$ is abelian). So it makes sense to associate to every component of
$D$ a cyclic subgroup
$H$ of $G$ and a generator $\psi$ of $H^*$. \\

For every pair of characters $\chi, \chi' \in G^*$, for every $H \in
\mathfrak{S}$ and for every $\psi \in S_H$ we can write
$$
\chi | H = \psi^{i_{\chi}}, ~~~\chi' | H = \psi^{i_{\chi'}}, ~~~
i_{\chi}, ~i_{\chi'} \in \{0, \ldots,
m_H - 1\},
$$
where $m_H$ is the order of $H$. One sets $\epsilon^{H,\psi}_{\chi,
\chi'} := 0$, if $i_{\chi} +
i_{\chi'} < m_H$ and $=1$ otherwise.

\begin{definition} Let $\pi : X \longrightarrow Y$ be an abelian
cover with group $G$. Moreover, assume
that  $X$ is normal and $Y$ smooth. Then the sheaves $L_{\chi}$,
$\chi \in G^*$, and the divisors
$D_{H, \psi}$ are called the {\em building data} of the covering.
\end{definition}

Now we are ready to formulate the following result of R. Pardini (cf.
\cite{Pardini1}, theorem $2.1$).

\begin{theo} {\bf(Pardini)} Let $G$ be a finite abelian group. \\ 1) 
Let $\pi : X
\longrightarrow Y$ be a covering with group $G$, where $X$ is a
normal variety and $Y$ is smooth and
complete. Then the building data of $\pi$ satisfy the following
linear equivalences:
$$ L_{\chi} + L_{\chi'} \equiv L_{\chi \chi'} + \sum_{H \in \mathfrak{S}}
\sum_{\psi \in S_H} \epsilon^{H,\psi}_{\chi, \chi'} D_{H, \psi_H}.
$$
2) Given any set of data $L_{\chi}$, $D_{H,\psi}$ satisfying the
above linear equivalences, there is a
unique (up to isomorphisms of Galois covers)  abelian cover \\
$\pi : X
\longrightarrow Y$ such that
$L_{\chi}$,
$D_{H,\psi}$ are its building data, if $X$ is normal.
\end{theo}

Sometimes, however, we shall also specify an Abelian cover through the
normalization of a singular Abelian cover. 

For instance, let
$P_1,
\ldots, P_5$ be five pairwise different points in
$\mathbb{P}^1$ and let $L$ be a divisor on
$\mathbb{P}^1$ such that
$$ 3L \equiv P_1 + P_2 + P_3 + P_4 + 2P_5 \equiv
\mathcal{O}_{\mathbb{P}^1}(6).
$$
Then there is a unique $\mathbb{Z}/3-$ covering $\pi : C_2 \longrightarrow
\mathbb{P}^1$ with branch locus $P_1 + P_2 + P_3 + P_4 + 2P_5$. We
remark that by Hurwitz' formula
the genus of $C_2$ has to be $3$. \\
    We calculate now the decomposition of $H^0(C_2,
\Omega^1_{C_2})$ into eigenspaces according
to the characters of $G$. We denote the characters of $G =
\mathbb{Z}/3$ simply by $0,1,2$. \\
Using the above theorem we can now calculate $L_{\chi}$. First of all, $L_0 =
\mathcal{O}_{\mathbb{P}^1}$ and $L_1 =
\mathcal{O}_{\mathbb{P}^1}(2)$. Using the above formula we
obtain $L_2 = L_1 + L_1 - P_5 =
\mathcal{O}_{\mathbb{P}^1}(3)$. By \cite{Pardini1}, prop. $4.1$, we know that
$$ (\pi_* \Omega^1_{C_2})^{\chi} = \Omega^1_{\mathbb{P}^1} \otimes
L_{\chi^{-1}}.
$$
Therefore we conclude that
$$ H^0(\Omega^1_{C_2}) = H^0(\Omega^1_{C_2})^0 \oplus H^0(\Omega^1_{C_2})^1
\oplus H^0(\Omega^1_{C_2})^2  =
$$
$$ H^0(\Omega^1_{\mathbb{P}^1}) \oplus H^0(\Omega^1_{\mathbb{P}^1} \otimes
\mathcal{O}(3)) \oplus H^0(\Omega^1_{\mathbb{P}^1} \otimes \mathcal{O}(2)) =
$$
$$ H^0(\mathbb{P}^1, \mathcal{O}(1)) \oplus H^0(\mathbb{P}^1,\mathcal{O}).
$$

\begin{theo} \label{counterexample} Let $C_1'$ be a general curve of genus
$3$ and let $C_1$ be a connected unramified covering of degree
$3$. Moreover, let $C_2$ be as
above. Then the following holds: \\ 1) $G :=
\mathbb{Z}/3$ operates freely on $C_1 \times C_2$; \\ 2)
infinitesimal Torelli fails for $S := C_1
\times C_2/ G$; \\ 3) $K_S$   is
quasi very ample, i.e., it gives a birational morphism which is 
a local embedding on the complement of a finite set.
\end{theo}

{\em Proof.} \\ We observe that, by proposition 2.3, we only have to 
show part $3)$.

  For
this we remark that $H^0(C_1,
\Omega^1_{C_1})^0$ has dimension $3$ while $h^0(C_1,
\Omega^1_{C_1})^1 = h^0(C_1, \Omega^1_{C_1})^2 = 2$.

Let $\eta$ be a $3$ - torsion element in
$Pic(C_1')$ corresponding to the unramified covering $C_1 \ra C'_1$.

Let $s_1, s_2
\in H^0(C_1, \Omega^1_{C_1})^1 =
H^0(C_1', K_{C_1'} + \eta)$ and $t_1, t_2 \in H^0(C_1,
\Omega^1_{C_1})^2 = H^0(C_1', K_{C_1'} + 2\eta)$
yield respective bases of these vector spaces.

  Let furthermore
$\sigma_1, \sigma_2 \in H^0(\Omega^1_{C_2})^1$,
resp. $ \sigma_3 \in H^0(\Omega^1_{C_2})^2$ yield respective bases of 
these vector spaces.
It is easy to see that the divisor of zeroes of $\sigma_3$ is
  the sum $ Q_1 + Q_2 + Q_3 + Q_4$, where $Q_i$ is the point lying over $P_i$.

Instead the linear series cut by $|\lambda_1 \sigma_1 + \lambda_2 
\sigma_2 = 0|$
has $Q_5$ as base point, and then yields residually the $g^1_3$ which gives the
triple covering of $\PP^1$.
It follows immediately that the curve $C_2$ is not hyperelliptic, whence
its canonical system is very ample. Observe then that
$$\{  s_1 \otimes \sigma_3, ~s_2 \otimes \sigma_3, ~t_1 \otimes \sigma_1, ~t_2
\otimes \sigma_2, ~t_1 \otimes \sigma_2, ~t_2 \otimes \sigma_1 \}
$$
is a basis of $H^0(S, \mathcal{O}(K_S))$. \\

We will use the following

\begin{lemma} \label{hilfslemma} Let $C_1'$ be a generic, (in
particular, non hyperelliptic) curve of genus $3$ and
let $\eta$ be a non trivial $3$-  torsion element of $Pic(C_1')$. Then 

1) the
linear system $|K_{C_1'} + \eta|$ is
base point free,

2) the morphism $ f : C_1' \ra \PP^1 \times \PP^1$ induced by the product
of the two pencils $|K_{C_1'} + \eta|$ and $|K_{C_1'} +  2 \eta|$ is
birational onto its image, which is a curve $\Gamma$ of bidegree $(4,4)$
having exactly 
$6$ ordinary double points as singularities, on the set $\M :=$
$\{ (x,y) | \  x \neq y, \ x,y \in \{
0,1, \infty \} \}$. 
\end{lemma}

Let's  postpone the proof of the lemma and infer from it the quasi 
very ampleness of the
canonical system of $S$.
\\
We observe first of all that the canonical system $K_S$ is base point free.

In fact, we can apply Lemma \ref{hilfslemma} twice, for $\eta$ and for $2 
\eta$, so that, given a point $(x,y) \in C_1 \times C_2$ we
may assume (after a change of basis) $s_2, t_2$ not to vanish on the image
of
$x$ in
$C'_1$.  And since the quotient of $C_2$ is rational, 
$H^0(C_2,
\Omega^1_{C_2})^0 = 0$ thus, the canonical system of a curve being 
base point free,  for
each choice of
$y$ there is a section
$\sigma_j$, for some $ 1 \leq j \leq3$, which does not vanish on $y$.

Let us now take two points  $(x,y), (x'',y'') \in C_1 \times C_2$
representing two distinct points in $S$, and let us assume that they are
not separated by the canonical map of $S$.

  If  $x= x''$, then the two points are separated since $\sigma_1,
\sigma_2, \sigma_3 $ yield the canonical map of $C_2$, which is non
hyperelliptic.

If $y = Q_5$ and the two points are not separated, it must be also
$y'' = Q_5$, since the last $4$ coordinates are then equal to zero.

But then $ (s_1(x) : s_2(x) ) = (s_1(x'') : s_2(x'') )$, which means that
$x'' $ and $ x$  map to  respective points $z'' $ and $ z$ such that
$z'' + z$ is contained in a divisor of $|K_{C_1'} + \eta|$. Since we have
two distinct points of $S$, $z'' \neq z$. We get therefore a singular
line in the canonical image $\phi_K(S)$, which is a quadruple line.

Assume now $y = Q_j , 1 \leq j \leq 4$ : then it must also be 
$y ''= Q_i , 1 \leq i \leq 4$ since the first two coordinates are zero.
If however $i \neq j$, then we may assume without loss of generality
$\sigma_1 (Q_j) = 0, \sigma_1 (Q_i) \neq 0$, a contradiction.

It follows that $y ''= y = Q_j , 1 \leq j \leq 4$, and we conclude as in
the previous case that $x'' $ and $ x$  map to  respective points 
$z'' \neq z$ such that
$z'',  z$ have the same image point under $|K_{C_1'} + 2 \ \eta|$. 
 We get thus $4$ more quadruple
lines in the canonical image $\phi_K(S)$. 

We  assume  now that $y '', y $ are no ramification points
($ y '',  y \neq Q_j , 1 \leq j \leq 5$). We may assume also as before
that $\sigma_1 (y) \neq 0, \sigma_2 (y) = 0$. We infer then immediately
that $y ''$ is in the $G$-orbit of $y$.

Thus, without loss of generality, and by the first step, we may assume 
 $y ''= y$ and that the respective image points $z, z''$ of $x, x''$ in
$C_1$ are distinct.

Since $\sigma_1 (y) \neq 0, \sigma_3 (y) \neq 0$, it follows that it must
be $f(z) = f(z'')$.
Fix then one of the $54$  pairs of points  $x, x''$  with the property
that $z \neq z''$, but $f(z) = f(z'')$.

We ask whether for  $y \in C_2$ the points $(x,y), ( x'', y)$ 
are not separated. Again after a change of basis, since   $f(z) = f(z'')$,
we may assume $s_1(x) = s_1(x'') = t_1(x) = t_1(x'') = 0$.
Then the question is whether $(s_2(x): t_2(x)) = (s_2(x''): t_2(x''))$.
In other words, the question we must answer is whether the subseries of
the canonical system of $C_1$, given by the $4$ sections $s_1,s_2, t_1,
t_2$ separates the two points $x, x''$.

We use here the fact that the point $f(z) = f(z'')$ is an ordinary double
point, i.e., with distinct tangents. 
There are local coordinates $(u,v)$ on $\PP^1 \times \PP^1$ such that
$f(z) = f(z'')$ corresponds to the origin. 
We show in the proof of lemma \ref{hilfslemma} that locally the triple
cover is obtained by extracting the third root  of $u/v$, $w^3= u/v$. 

Accordingly, the $4$ sections $s_1,s_2, t_1, t_2$ are locally expressed by
$1, u, w, wv$. If $w(x) = w(x'')$, then we would have $(u/v)(z) =
(u/v)(z'')$, a contradiction.

We have thus proven that for $C'_1$ generic the canonical system of $S$
yields a morphism which is injective except for the curves images of
$C_1 \times \{ Q_j\} \ , 1 \leq j \leq 5  $.

There remains to prove that the canonical map is a local 
embedding except
at finitely many points.

Lemma \ref{hilfslemma} shows that for a general curve $C$ of genus $3$,
the mapping $f : C \ra \P^1 \times \PP^1$ ,
corresponding to the sum of the two linear systems
$|K_{C} + \eta|$ and $|K_{C} + 2 \eta|$, is everywhere a local embedding.

We will prove that the canonical system $|K_S|$ separates
tangent vectors at a point  corresponding to $(x,y) \in C_1 \times C_2$,
unless $y$ is  a ramification point for the triple cover of $\PP^1$,
 and $x$
is  a ramification point for one of the two $g^1_4$'s $|K_{C} + 
\eta|$ and $|K_{C} + 2 \eta|$.

Without loss of generality we may then take a point $(x,y) \in C_1 
\times C_2$ and we assume 
 $s_1, t_1$ to vanish on  $x$ and $s_2, t_2$ not to vanish on $x$.

We look for two sections yielding two  curves which are smooth at 
$(x,y)$ and have
distinct tangents.

If we have that $s_1$ vanishes of order exactly $1$ at $x$ and 
$\sigma_3 (y) \neq 0$,
we are done, since there is always a section $\sigma_i$ vanishing 
simply at $y$,
and $s_1 \sigma_3 $, $s_2 \sigma_i$ give two curves with vertical,
 respectively
horizontal tangent.

Similarly, if $t_1$  vanishes of order exactly $1$ at $x$, and $y \neq 
Q_5$, then
we may assume  $\sigma_1 (y) \neq 0$, and   $\sigma_2$ to vanish 
simply on $y$, unless we
are in a branch point $Q_i$. For $ i \leq 4$, however, $\sigma_3$ 
vanishes simply and
we are therefore done.

\qed for Theorem \ref{counterexample}

{\em Proof of lemma \ref{hilfslemma}.} \\ The moduli space of curves
of genus $3$ has dimension $6$ and
the hyperelliptic curves form a five dimensional algebraic subset.
Hence we can suppose that $C_1'$ is
not hyperelliptic. Therefore $C_1'$ is canonically embedded as a
plane quartic in $\mathbb{P}^2$. Let
$\eta$ be an element of $Pic(C_1')_3$. We note that
$P$ is a base point of the linear system $|K_{C_1'} + \eta|$ if and only if
$$ H^0(C_1', \mathcal{O}(K + \eta)) = H^0(C_1', \mathcal{O}(K + \eta - P)).
$$
Since $dim H^0(C_1', \mathcal{O}(K + \eta)) = 2$ this is equivalent to
$$ dim H^1(C_1', \mathcal{O}(K + \eta - P)) = 1.
$$
Since $H^1(C_1', \mathcal{O}(K + \eta - P)) \cong H^0(C_1', \mathcal{O}(P -
\eta))^*$, it follows that there is a point $P'$ such that $P - \eta
\equiv P'$. Therefore $3P \equiv
3P'$. By Riemann - Roch we have
$$ dim H^0(C_1', \mathcal{O}(K - 3P)) =
$$
$$ deg(K - 3P) + 1 - g(C_1') + dim H^1(C_1', \mathcal{O}(K - 3P) = 1
+ 1 - 3 + 2 = 1,
$$
and there is a point $Q$ such that $Q \equiv K - 3P \equiv K - 3P'$. \\

Geometrically this means that, considering $C_1'$ as a plane quartic,
$C_1'$ has two inflection points $P$, $P'$, such that the tangent
lines to these points intersect in $Q
\in C_1'$. \\
Let now $p$, $q$, $p'$ be three non collinear points in
$\mathbb{P}^2$. Then the dimension of the group of
automorphisms of $\mathbb{P}^2$ leaving the three points fixed has
 dimension $2$. The
quartics in $\mathbb{P}^2$ form a linear system of dimension $14$.
Imposing that a plane quartic
contains the point $q$ is one linear condition. Moreover, the
condition that the line containing $p$
and $q$ is a tritangent to the quartic gives further three linear
conditions as well as the condition
that the line containing $p'$ and $q$ is a tritangent to the quartic.
Therefore the linear subsystem of
quartics $C$ having two inflection points $P$, $P'$, such that the
tangent lines to these points
intersect in $Q \in C$ has dimension $14 - 2 - 3 - 3 - 1 = 5$, whence
they are special.

Part 1) is therefore proven and we set for convenience of notation $C: =
C_1'$. By 1) we may assume that both linear systems $|K_C + \eta|$
and $|K_C + 2 \ \eta|$ are base point free.

Even if not strictly needed for our purposes, we try to describe how one
reaches the conclusion that the curve $C$ must be birational to  a
curve $\Gamma \subset\PP^1 \times \PP^1$ enjoying the properties stated in
2).

More generally, for $L \in Pic^0(C)$, $L$ general, we obtain a morphism
$f_L : C \ra \PP^1 \times \PP^1$ corresponding to the linear system sum of 
 $|K_C + L|$ and of  $|K_C -  L|$.

Set $\Gamma : = f_L(C)$. Then either $\Gamma$ is a curve of bidegree
$(4,4)$, or $ deg (f_L) = 2$ and $\Gamma$ is a curve of bidegree
$(4,4)$, since if $ deg (f_L) = 4$ then $ L \equiv -L$ and $L$ is thus
of $2$-torsion.

We will assume that $f_L$ is birational, else $C$ is either hyperelliptic,
or a double cover of an elliptic curve, which is then special because it
is branched in $4$ points, and we get only a family of dimension $5$ in
the moduli space.

Let $P_1, \dots P_m$ be the (possibly infinitely near) singular points of
$\Gamma$.

Then, $H_1$, $H_2$ being the respective divisors of a vertical, and of a
horizontal line in $\PP^1 \times \PP^1$,
$\Gamma \in | 4 H_1 + 4 H_2 - \Sigma_{i=1, \dots m} r_i P_i |$.

By adjunction, the canonical system of $C$ is cut by the series 
$| 2 H_1 + 2 H_2 - \Sigma_{i=1, \dots m} (r_i - 1) P_i |$.

Since $C$ has genus $3$, we obtain 
$\Sigma_{i=1, \dots m} r_i (r_i - 1) = 12$. Whence, either $\Gamma$ has
just an ordinary $4$-uple point or all the multiplicities are at most $3$.

In the former case, however, we take local coordinates $(u,v)$ at the 
$4$-uple point and see that the canonical map is given by $(u^2 v : v^2 u
: u^2 v^2) = (u : v : uv)$ and the curve $C$ is hyperelliptic.

In the latter case, we observe that $\Gamma$ lies on a regular surface,
whence  we have that the bicanonical system is obtained as the restriction
of the linear system $| 4 H_1 + 4 H_2 - \Sigma_{i=1, \dots m} 2 (r_i - 1)
P_i |$ ( in fact on the blown up surface $X$, we have $-1$ divisors $E_i$
and  $ H^1 (\hol_X( - \Sigma_{i=1, \dots m} (2 - r_i )E_i )) = 0$
since $r_i \leq 3$, by Ramanujam's vanishing ).

We exploit at this point that $H_1$ pulls back to $K_C + L$, $H_2$ pulls
back to $K_C - L$, thus there is a subseries of the bicanonical series
cut by $ H_1 + H_2$. We infer the existence of a curve $G \in  |3 H_1 + 3
H_2  - \Sigma_{i=1, \dots m} 2 (r_i - 1) P_i |$.

Assuming that $G$ is reduced, we see that $G$ has geometric genus 
$ g(G) = 4 -  \Sigma_{i=1, \dots m} (2 r_i - 3) (r_i - 1) P_i  \leq -2$,
and indeed $ g(G)  \leq -5 $ if there is a triple point for $\Gamma$.
Yet, $G$ has at most $6$ components, thus $ g(G)  \geq -5 $, and if
equality holds, $G$ consists of $3$ vertical and $3$ horizontal lines, in
particular it has no points of multiplicity $4$.

It follows then that  all the $r_i = 2$, and $G$ has three components,
which are rational. 

{\bf General case. } $G = q_1 \cup q_2 \cup q_3$, where $q_j$ is of
bidegree $(1,1)$, the curves $q_i$ and $q_j$ meet transversally in 
$2$ points $P_{i.j}, P'_{i,j}$.

Fixing these three conics $q_j$, one takes a general curve $\Gamma \in
 | 4 H_1 + 4 H_2 - \Sigma_{i,j} 2 (P_{i.j} + P'_{i,j})|$. 

Let $X$ be the blow up of $\PP^1 \times \PP^1$ in these $6$ points, and let
$\tilde{G} ,
\tilde{\Gamma}$ be the proper transforms of $G$, resp. $\Gamma$. By the
exact sequence 
$$ 0 \ra \hol_X (H_1 + H_2)  \ra  \hol_X ( \tilde{\Gamma}) \ra 
 \hol_{\tilde{G}} ( \tilde{\Gamma}) \ra 0$$
we obtain  dim $ | \tilde{\Gamma}| = 4 + 3 -1 = 6$. Since we have moreover 
$ 3 + 3 + 3 $ moduli for the $3$ conics $ q_j$, after subtracting $6$
moduli for the automorphisms of $\PP^1 \times \PP^1$, we obtain a family
of dimension $9 = 6 + 3$, which is the expected dimension since $6$ is the
dimension of the moduli space of curves of genus $3$, and $ 3 = { \rm dim}
\ Pic^0(C)$.\\

{\bf 3 - Torsion case. } Assume now that $ L = \eta$ is a non trivial
$3$-torsion divisor: then $ 3 H_i \equiv 3 K_C$, whence we expect curves
$\Delta_1 \in | 6 H_1 + 3 H_2 - \Sigma_{i=1, \dots m} 3 (r_i - 1)
P_i |$ , $\Delta_2 \in | 3 H_1 + 6 H_2 - \Sigma_{i=1, \dots m} 3 (r_i - 1)
P_i |$. Intersecting on $X$ the proper transform of $\Delta_i$ with the
proper transform of each component $q_j$ of
$G$ we get intersection number $ 9 - 12 = -3$, thus $G$ should be
contained in $\Delta_i$.

But then $\Lambda_1 : = \Delta_1 - G$, $\Lambda_2 : =\Delta_2 - G$
provide divisors in $ | 3 H_1 - \Sigma_{i=1, \dots m}  (r_i - 1)
P_i |$, respectively in  $ | 3 H_2 - \Sigma_{i=1, \dots m}  (r_i - 1)
P_i |$. Again , the intersection number of the proper transforms of $G$
and $\Lambda_i$ are negative $ = 6 - 12 = -6$ , thus we conclude that
$ G = \Lambda_1  + \Lambda_2 $, and that each $q_j$ is reducible.

The curve $\Gamma$ has $6$ double points which lie on these $6$ lines,
but since the bidegree of  $\Gamma$ is $(4,4)$ there are at most $2$
singular points on each line, and without loss of generality, since we
assumed $G$ to be reduced, we may assume that the $6$ points are as in the
statement of the Lemma, namely, they are the points of the set $\M : =$
$\{ (x,y) | \  x \neq y, \ x,y \in \{ 0,1, \infty \} \}$.\\

{\bf Claim 1.} The general element $\Gamma \in | 4 H_1 + 4 H_2 - \M|$
is irreducible and has only ordinary double points as singularities.

 {\em Proof of claim 1}: The general element has only  double points as
singularities, since we may just consider the subsystem $G$ + 
$|H_1 + H_2|$. Since the selfintersection of the proper transforms on $X$
of the curves in the linear system equals $ 32 - 24= 8$, the system is not
composed of a pencil, and the general member is irreducible.

{\bf Claim 2.} For any curve $\Gamma \in | 4 H_1 + 4 H_2 - \M|$ we get an
unramified cyclic triple cover $Y$ of the normalization $C$  of $\Gamma$ by
taking as $Y$ the normalization of the inverse image of $\Gamma$ in the
cyclic triple cover $W$ of
$\PP^1 \times \PP^1$ branched on $G$, and such that $\Lambda_1$ is the
divisor $D_1$ in Pardini's notation, $\Lambda_2$ is the divisor $D_2$ .\\

{\em Proof of claim 2}: $G$ intersects $\Gamma$ only in the double
points. Once we pull back the triple covering to $C$, for each point $p$
lying over a double point, the point $p$ apears in the branch divisor with
multiplicity $ 1-1 = 0$. 

Thus $ Y \ra C$ is unramified.\\

{\bf Final Observation .} The linear system $| 4 H_1 + 4 H_2 - \M|$ has
projective dimension equal to $6$, as previously indicated. We get
therefore a $6$ dimensional family, which is therefore dominant onto the
moduli space of pairs $ (C, \eta)$ as above.

\qed

\begin{rem} The family of surfaces  considered in
Theorem \ref{counterexample} has dimension $8$.

We observe that the canonical image $\phi_K(S)$ has $S$ as its
normalization. The inverse images of the singular curves of $\phi_K(S)$ 
are exactly the curves $ C_1 \times \{ Q_j\}, 1 \leq j \leq 5 $.

Therefore, the geometry of the canonical image $\phi_K(S)$  recovers
$C_1$, which is the "missing" Hodge structure. 
\end{rem}

Replacing now the curve $C_2$ by a triple cover of
$\mathbb{P}^1$ ramified in more points than
in the previous example, we can produce an infinite series of
examples of surfaces having quasi very ample
canonical bundle for which the infinitesimal Torelli theorem fails.
The construction goes as follows.\\

Let $k$ be any natural number bigger or equal to two  and let  $P_1,
\ldots, P_{3k-2}, P_{3k-1}$ be pairwise different points in
$\mathbb{P}^1$. We consider the triple cover
$C_k$ of  $\mathbb{P}^1$ ramified along
$$ P_1 + \ldots + P_{3k-2} + 2P_{3k-1}.
$$
Then the genus of $C_k$ is equal to $3k - 3$. Using Pardini's
formulae we can calculate the
decomposition of $H^0(\Omega^1_{C_2})$ into eigenspaces according
to the characters of $G =
\mathbb{Z}/3$. We get: $L_0 =
\mathcal{O}_{\mathbb{P}^1}$, $L_1 = \mathcal{O}_{\mathbb{P}^1}(k)$ and $L_2 =
\mathcal{O}_{\mathbb{P}^1}(2k - 1)$ and therefore
$$ H^0(\Omega^1_{C_k}) = H^0(\Omega^1_{C_k})^0 \oplus H^0(\Omega^1_{C_k})^1
\oplus H^0(\Omega^1_{C_k})^2  =
$$
$$ H^0(\Omega^1_{\mathbb{P}^1}) \oplus H^0(\Omega^1_{\mathbb{P}^1} \otimes
\mathcal{O}(2k-1)) \oplus H^0(\Omega^1_{\mathbb{P}^1} \otimes
\mathcal{O}(k)) =
$$
$$ H^0(\mathbb{P}^1, \mathcal{O}(2k-3)) \oplus
H^0(\mathbb{P}^1,\mathcal{O}(k-2)).
$$
Using the same argument as in theorem \ref{counterexample} (note that
we only use that the curve $C_2$
is not hyperelliptic and that $C_2/G =
\mathbb{P}^1$) we obtain the following

\begin{theo} Let $k$ be any natural number bigger or equal to $2$ and let
$C_1'$ be a general (non hyperelliptic) curve of genus $3$ and let
$C_1$ be a connected unramified
covering of $C_1'$  of degree $3$. 
Moreover, let $C_k$ be as above. \\ Then
the following hold: \\ 1) $G :=
\mathbb{Z}/3$ operates freely on $C_1 \times C_k$; \\ 2)
infinitesimal Torelli fails for $S_k := C_1
\times C_k/ G$; \\ 3) $K_{S_k}$ is  quasi very ample.
\end{theo}

\section{Torelli fibres}

In this section we want to analyse the above series of examples in
the following sense. We know that the infinitesimal Torelli map
 fails to be injective for
reasons of dimension of the eigenspaces, but we want to calculate
 the dimension of the kernel of the infinitesimal Torelli
map.

We start by computing the respective dimensions of
$H^1(S_k,\Omega^1_{S_k})$, $ H^0(S_k,\Omega^2_{S_k})$ and  $
H^1(S_k,\Omega^1_{S_k} \otimes
\Omega^2_{S_k})$,  and end counting the dimension of the moduli space of
the surfaces $S_k$.\\

Let $k$ be any natural number bigger or equal to two  and let  $P_1,
\ldots, P_{3k-2}, P_{3k-1}$ be pairwise different points in
$\mathbb{P}^1$. We consider the triple cover
$C_k$ of  $\mathbb{P}^1$ ramified along
$$ P_1 + \ldots + P_{3k-2} + 2P_{3k-1}.
$$
Then we define: $S_k := (C_1 \times C_k)/G$, where $G :=
\mathbb{Z}/3\mathbb{Z}$. We prove the following

\begin{prop} \label{dimensions} Under the above hypothesis we have: \\
1) $h^0(S_k,\Omega^2_{S_k}) = 6k - 6$, \\
2) $h^1(S_k,\Omega^1_{S_k}) = 12k - 10$, \\
3) $h^1(S_k,\Omega^1_{S_k} \otimes \Omega^2_{S_k}) = h^1(S_k,
T_{S_k}) = 3k + 2$.

\end{prop}

\begin{rem} In particular, it follows from the above proposition that the
infinitesimal Torelli theorem
does not fail for reasons of dimension of the vector spaces in question,
but only because of the imposed symmetry.
\end{rem}

{\em Proof.} For the dimensions of the eigenspaces of
$H^0(\Omega^1_{C_1})$, resp. $H^0(\Omega^1_{C_k})$,
we have the following: \\
a) $H^0(C_1, \Omega^1_{C_1})^0 = \mathbb{C}^3$, \\
b) $H^0(C_1, \Omega^1_{C_1})^1 = \mathbb{C}^2$, \\
c) $H^0(C_1, \Omega^1_{C_1})^2 = \mathbb{C}^2$, \\
d) $H^0(C_k, \Omega^1_{C_k})^0 = 0$, \\
e) $H^0(C_k, \Omega^1_{C_k})^1 = H^0(\mathbb{P}^1, \mathcal{O}(2k-3)) =
\mathbb{C}^{2k-2}$, \\
f) $H^0(C_k, \Omega^1_{C_k})^2 = H^0(\mathbb{P}^1, \mathcal{O}(k-2)) =
\mathbb{C}^{k-1}$. \\

Using the decomposition of the cohomology groups of $Z_k := C_1
\times C_k$ into eigenspaces
corresponding to the characters of $G$ we obtain
$$ (H^0(C_1, \Omega^1_{C_1})^1 \otimes H^0(C_k, \Omega^1_{C_k})^2)
\oplus (H^0(C_1, \Omega^1_{C_1})^2
\otimes H^0(C_k, \Omega^1_{C_k})^1) =
\mathbb{C}^{6k-6},
$$
which proves 1). \\

Moreover,
$$ ( \ H^1(C_1,\Omega^1_{C_1}) \otimes H^0(C_k, \mathcal{O}_{C_k}) \
) \oplus ( \ H^0(C_1,\mathcal{O}_{C_1})
\otimes H^1(C_k, \Omega^1_{C_k}) \ ) \oplus
$$
$$
\bigoplus_{\chi \in G^*} (H^0(C_1,\Omega^1_{C_1})^{\chi}  \otimes H^1(C_k,
\mathcal{O}_{C_k})^{\chi^*}) \oplus \bigoplus_{\chi \in G^*}
(H^1(C_1,\mathcal{O}_{C_1})^{\chi} \otimes H^0(C_k,
\Omega^1_{C_k})^{\chi^*}) =
$$
$$
\mathbb{C}^2 \oplus (H^0(C_1, \Omega^1_{C_1})^1 \otimes H^1(C_k,
\mathcal{O}_{C_k})^2) \oplus (H^0(C_1, \Omega^1_{C_1})^2 \otimes H^1(C_k,
\mathcal{O}_{C_k})^1) \oplus
$$
$$
\oplus (H^1(C_1, \mathcal{O}_{C_1})^1 \otimes H^0(C_k, \Omega^1_{C_k})^2)
\oplus (H^1(C_1, \mathcal{O}_{C_1})^2 \otimes H^0(C_k, \Omega^1_{C_k})^1) =
$$
$$
\mathbb{C}^2 \oplus \mathbb{C}^{4k-4} \oplus \mathbb{C}^{2k-2} \oplus
\mathbb{C}^{2k-2} \oplus \mathbb{C}^{4k-4} = \mathbb{C}^{12k-10},
$$

and this proves 2). \\

Finally, we see that
$$ 
H^1(S_k, \Omega^1_{S_k} \otimes \Omega^2_{S_k}) =
$$
$$
\bigoplus_{\chi \in G^*} (H^0(C_1, \mathcal{O}_{C_1}(2K_{C_1}))^{\chi}
\otimes H^1(C_k, \mathcal{O}_{C_k}(K_{C_k}))^{\chi^*}) \oplus
$$
$$
\bigoplus_{\chi \in G^*} (H^1(C_1, \mathcal{O}_{C_1}(K_{C_1}))^{\chi}
\otimes H^0(C_k, \mathcal{O}_{C_k}(2K_{C_k}))^{\chi^*}).
$$
Since $H^1(C_1, \mathcal{O}_{C_1}(K_{C_1}))$ and $ H^1(C_k,
\mathcal{O}_{C_k}(K_{C_k}))$ are automatically invariant under $G$, we get that
$$ 
H^1(S_k, \Omega^1_{S_k} \otimes \Omega^2_{S_k}) = (H^0(C_1,
\mathcal{O}_{C_1}(2K_{C_1})^G \otimes  H^1(C_k, \mathcal{O}_{C_k}(K_{C_k})))
\oplus
$$
$$
\oplus (H^1(C_1, \mathcal{O}_{C_1}(K_{C_1}) \otimes H^0(C_k,
\mathcal{O}_{C_k}(2K_{C_k}))^G).
$$
Using
$$ 
H^0(C_1, \mathcal{O}_{C_1}(2K_{C_1}))^G = H^0(C'_1,
\mathcal{O}_{C'_1}(2K_{C'_1})) = \mathbb{C}^6
$$
and
$$ h^1(C_1, \mathcal{O}_{C_1}(K_{C_1})) = h^1(C_k,
\mathcal{O}_{C_k}(K_{C_k})) = 1,
$$
we see that we have proven 3) as soon as we have shown that
$$ 
H^0(C_k, \mathcal{O}_{C_k}(2K_{C_k}))^G = \mathbb{C}^{3k-4}.
$$
This will be done in the following lemma.

\qed \\

First we recall the definition of almost simple cyclic coverings (cf.
\cite{cat89}, page 309). Let
$Y$ be an algebraic manifold and let
$\mathcal{L} = \mathcal{O}_Y(F)$ be the invertible sheaf which is the
sheaf of sections of a line bundle $\LL$ on
$Y$. Assume that there are given reduced effective
divisors $\Delta_0$, $\Delta_{\infty}$ on $Y$, which are disjoint and
it holds $\Delta_0 \equiv
\Delta_{\infty} + nF$.

\begin{definition} The {\em almost simple cyclic cover} associated to
$(Y,\mathcal{L},\Delta_0,\Delta_{\infty})$ is the subvariety $X
\subset \mathbb{P} (\LL \oplus
\mathbb{C}_Y)$ defined by the equation $z_1^n \delta_{\infty} =
\delta_0 z_0^n$, where $\delta_{\infty}$, $\delta_0$ are sections defining
$\Delta_{\infty}$, resp. $\Delta_0$, and $z_1$, $z_0$ are respective
linear coordinates on the fibres of
$\LL$, resp. the trivial line bundle $\mathbb{C}_Y$.
\end{definition}

We have the following :

\begin{lemma} Let $\pi :X \longrightarrow Y$ be an
almost simple cyclic covering of
degree $n$ with smooth branch divisors $ \Delta_0 , \Delta_{\infty}$ and
set
$G := \mathbb{Z} / n$. Then:
$$ 
H^0(X, \mathcal{O}_{X}(2K_X))^G = H^0(Y,
\mathcal{O}_Y
(2K_Y + \Delta_{\infty} + \Delta_0)).
$$
\end{lemma}

{\em Proof.} \\
$X$ sits in the projective bundle $\mathbb{P}: = \mathbb{P} (\LL \oplus
\mathbb{C}_Y) = Proj (\mathcal{L}^{-1} \oplus \hol_Y)$
and is linearly equivalent to $ n H + p^* \Delta_0$, where $H$ is the
hyperplane divisor $H = div (z_0)$.

The canonical divisor of $\PP$ equals, by the relative Euler sequence,
$ - 2 H + p^* ( K_Y - F)$.
 Thus, by adjunction, $  K_X$ is the
restriction of the divisor 
$$ (n-2) H + p^* ( K_Y - F + \Delta_0).$$

If we set $ E_i : = div (z_i)$,  we  may write  
$$  2 K_X \equiv (n-2) E_0 + (n-2) E_1 +  p^* ( 2 K_Y + \Delta_0 +
\Delta_{\infty}).$$

It suffices to show that each invariant section of  
$H^0(X, \mathcal{O}_{X}(2K_X))$ vanishes on $E_i$ of multiplicity $n-2$
and after being divided by $ (z_0 z_1)^{n-2}$ yields an invariant section
of $ H^0(X,
\mathcal{O}_X
(p^*(2K_Y + \Delta_{\infty} + \Delta_0))).$

This follows by a local calculation, since, if we have that
 $y_0 = 0$ is a local equation of $\Delta_0$, we have that $z_1/ z_0 : =
x_0$ is a local equation for the ramification locus $E_1$, 
and we can complete $y_0 $  to local coordinates $(y_0 , y_1, \dots )$
on $Y$.

Note that  $\pi$ is locally given by an equation $x_0^n = y_0$,
hence $(x_0 , y_1, \dots )$ are local coordinates on $X$.

A differential form
$f(x_0, y_1, \dots )(dx_0 \wedge dy_1 \wedge \dots)^2$ is $G$ - invariant
if and only if
$f(x_0, y_1, \dots) = x_0^{n-2}
\gamma (y_0 , y_1, \dots )$. 

Hence
$$ f(x_0, y_1, \dots )(dx_0 \wedge dy_1 \wedge \dots)^2 = \gamma (y_0 ,
y_1, \dots) x_0^{n-2} (dx_0 \wedge dy_1 \wedge \dots)^2 =$$
$$ = (n)^{-2} \gamma (y_0 ,
y_1, \dots) x_0^{n-2}  (dy_0 \wedge dy_1 \wedge \dots)^2
\frac{1}{x_0^{2n-2}} = (n)^{-2}   (dy_0 \wedge dy_1 \wedge \dots)^2
\cdot \frac{\gamma (y_0 , y_1, \dots) }{y_0} .
$$
The same calculation holds for $ \Delta_{\infty}$ and we are done.

\qed \\

In our particular case this implies immediately:

\begin{cor} \label{2kinv} Let $k$ be any natural number bigger or
equal to two  and let  $P_1, \ldots,
P_{3k-2}, P_{3k-1}$ be pairwise different points in $\mathbb{P}^1$. We
consider the triple cover $C_k$ of
$\mathbb{P}^1$ ramified along
$$ 
P_1 + \ldots + P_{3k-2} + 2P_{3k-1}.
$$
Then we have
$$ 
H^0(C_k, \mathcal{O}_{C_k}(2K_{C_k}))^G \cong \mathbb{C}^{3k-4}.
$$
\end{cor}

{\em Proof.} \\ 
Here $\Delta_0 = P_1 + \ldots P_{3k-2}$ and
$\Delta_{\infty} = P_{3k-1}$ and we get
$$ 
H^0(C_k, \mathcal{O}_{C_k}(2K_{C_k}))^G = H^0(\mathbb{P}^1,
\mathcal{O}_{\mathbb{P}^1} (2K + \Delta_0 + \Delta_{\infty})).
$$
$$ 
H^0(\mathbb{P}^1, \mathcal{O}_{\mathbb{P}^1} (- 4 + 3k - 1)) =
H^0(\mathbb{P}^1,
\mathcal{O}_{\mathbb{P}^1} (3k - 5)) = \mathbb{C}^{3k-4}.
$$
   \qed \\

We calculate now the dimension of the fibres of the
infinitesimal Torelli map. 

\begin{lemma} Let $k \geq 2$ be a natural number and consider the surfaces
$S_k$ as above. Then the fibres of the period map $\Phi_2$ have (at
each point)  dimension at least
two.
\end{lemma}

{\em Proof.} We know that
$$ dim \Phi(\bigoplus_{\chi \in G^* : H^0(C_k,\Omega^1_{C_k})^{\chi} \neq 0}
H^0(C_1, \Omega^1_{C_1})^{\chi}
\otimes H^0(C_1, \Omega^1_{C_1})^{\chi^*})
\leq
$$
$$ dim (H^0(C_1, \Omega^1_{C_1})^1 \otimes H^0(C_1, \Omega^1_{C_1})^2) = 4,
$$
whereas $dim H^0(C_1, \Omega^1_{C_1} \otimes \Omega^1_{C_1}) = 6$.
By the structure of the
infinitesimal Torelli map (cf. proof of theorem
\ref{thm1}) our claim is proven.
\qed \\

\begin{rem} Let $k \geq 2$ be a natural number and let
$\mathcal{S}_k$ be the family of surfaces we
have constructed above, which give  an irreducible connected component
of the moduli space of surfaces of general type by Theorem C of
\cite{cat00}.

Observe that the number of moduli of $S_k$ is equal to the sum of the
number of moduli of $C_1$ ($C_1$ is a curve of genus three,
hence has $6$ moduli) and of the number of moduli of the
triple cover $C_k \ra \mathbb{P}^1$ (it is ramified in $3k-1$ points,
whence it has $3k-1-3 = 3k-4$ moduli). From
this we conclude that $dim \mathcal{S}_k = 3k + 2$. \\ In particular,
since by proposition
\ref{dimensions}, 3) $h^1(S_k, T_{S_k}) = 3k + 2$, the family
$\mathcal{S}_k$ yields a generically smooth moduli space, and a smooth
base of the Kuranishi family.
\end{rem}
    
We can summarize our results now in the following theorem.

\begin{theo} For any natural number $k \geq 2$ there exists a $3k+2$-
dimensional family $\mathcal{S}_k$ of surfaces 
 such that the following holds:\\ 
1) for each
$S_k \in \mathcal{S}_k$ the
infinitesimal Torelli map
$$ 
d\Phi_2 : H^1(S_k,T_{S_k}) \longrightarrow Hom(H^{2,0}(S_k), H^{1,1}(S_k))
$$ 
has a kernel of dimension at least two, i.e., the period map
$\Phi_2$ has all  fibres of  dimension strictly positive and $\geq
2$. \\  2) for a general $S_k
\in \mathcal{S}_k$ the canonical divisor $K_{S_k}$ is 
quasi very ample;
\\ 
3) $h^1(S_k,T_{S_k}) = 3k+2$; in particular, $\mathcal{S}_k$ yields a
generically smooth irreducible connected
component of the moduli space.
\end{theo}

\begin{rem} 
We observe that, since our examples are irregular
surfaces, there is still another period
map $\Phi_1$ associated to the weight one Hodge structure on $S_k$.
The differential of this map is
given by
$$ 
d\Phi_1 : H^1(S_k,T_{S_k}) \longrightarrow Hom(H^{1,0}(S_k), H^{0,1}(S_k)),
$$
and its injectivity is equivalent to the surjectivity of
$$ 
d\Phi_1^* : H^0(S,\Omega^1_S) \otimes H^1(S,\Omega^2_S)
\longrightarrow H^1(S,\Omega^1_S \otimes
\Omega^2_S).
$$
Using again the explicit description of $S_k$ as in the proof of theorem
\ref{thm1}, we see that $d\Phi_1^*$ is surjective if and only if the two maps
$$
\varphi_i : H^0(C_i, \Omega^1_{C_i})^G \otimes H^0(C_i, \Omega^1_{C_i})^G
\longrightarrow H^0(C_i,\Omega^1_{C_i} \otimes \Omega^1_{C_i})^G
$$
are surjective. Since $C_2^k$ is a triple cover of $\mathbb{P}^1$ the
above map for $i=2$ is obviously
not surjective. On the other hand $\varphi_1$ is given by the natural map
$$  
H^0(C'_1, \Omega^1_{C'_1}) \otimes H^0(C'_1, \Omega^1_{C'_1})
\longrightarrow H^0(C'_1,\Omega^1_{C'_1} \otimes \Omega^1_{C'_1})
$$
and this map is surjective if and only if $C'_1$ is not hyperelliptic. \\
In the case that $C'_1$ is hyperelliptic, the cokernel of $\varphi_1$
has dimension $1$ (in general
$g-2$) and is transversal to the cokernel of
$\Phi_2$.
\end{rem}

\section{Global double Torelli for surfaces isogenous to a product}

Recall that the global Torelli theorem for algebraic curves says that 
the isomorphism
class of a curve $C$ is completely determined by the isomorphism 
class of the datum of
the integral cohomology algebra $H^* (C, \mathbb{Z})$ together with the 
Hodge decomposition
$$
H^1 (C, \mathbb{Z}) \otimes_{\mathbb{Z}} \CC =  H^{1,0} (C) \ \oplus H^{0,1} (C),
$$
where $H^{1,0} (C) = H^0(C, \Omega^1 _{C} )$ and $ H^{0,1} (C) = 
\overline{H^{1,0}
(C)}.$

We shall shortly say that $C$ is determined by its integral Hodge structure.

A fortiori, $C$ is determined by its topological type given together 
with its  Hodge
decomposition.

In this section we want to address the question: when does global
double Torelli hold for $S = C_1
\times C_2 / G$ ? ``Double'' means that the datum of the respective 
Hodge structures
of  weight one and two of $S$, together with the oriented topological type,
should determine the isomorphism class of $S$.

To be more precise and explicit, we first observe the following:

\begin{lemma} Let $G$ be a finite group acting on a Riemann surface
$C$ of genus $g \geq 2$. Then the pair consisting of the orbifold 
group exact sequence and of the Hodge structure of $C$
$$ 
(1 \rightarrow \pi_1(C) \rightarrow \pi_1^{orb} \rightarrow G
\rightarrow 1 ,
~H^0(\Omega^1_C) \subset H^1 (C, \mathbb{C}) = H^1 (\pi_1(C), \mathbb{C}))
$$
determines the holomorphic action of $G$ on $C$.
\end{lemma}

{\em Proof.}

Inner conjugation on  $\pi_1(C)$ of a lift of an element  $\gamma \in G$  
provides an inclusion of $G$ in the mapping class group
$G
\ra Out (\pi_1(C)) = Map_g$ of
$C$, which, together with the natural orientation provided by the half line
$H^0(\Omega^1_C)) \wedge \overline{~H^0(\Omega^1_C))} \subset H^1 (C,\mathbb{R})$,
determine the oriented topological type of the action. The Hodge 
structure determines the
complex structure on the curve which makes the $G$-action holomorphic.

\qed

Assume now that $S$ is a surface isogenous to a product and not of 
mixed type: by
\cite{cat00}, p. 25/26 it follows  that $\pi_1(S)$ determines the exact
sequence
$$ 1 \rightarrow \Pi_{g_1} \times \Pi_{g_2} \rightarrow \pi_1(S)
\rightarrow G \rightarrow 1
$$
which in turn determines the two orbifold $\pi_1$- exact sequences
$$ 1 \rightarrow \Pi_{g_i} \rightarrow \pi_1^{orb}(C_i' - B, m''_i)
\rightarrow G \rightarrow 1.
$$

Therefore the double Torelli question reduces to the following: \\ 
Does the
Hodge structure on $S$ determine the Hodge
structures on $C_1$ and $C_2$? 

\begin{rem} Let $\Delta : G \rightarrow G \times G$ be 
the diagonal inclusion.

Then, since $G$ is Abelian,  $G = G \times G/ \Delta (G)$
acts on the cohomology algebra $H^*(S)$ and we get the following
decompositions according to the
characters  (here, $\chi^*$ denotes the character of the dual 
representation, i.e.,
$\chi^*$ is the inverse of $\chi$, and $ \chi^* = (\chi)^{-1} = 
\overline{\chi}$):

$ H^1(S, \mathbb{Z}) $, which is a free $\ZZ$- module, contains as a
submodule of finite index 

$ H^1(C_1 \times C_2, \mathbb{Z})^G = H^1(C_1, \mathbb{Z})^G\oplus H^1(C_2,
\mathbb{Z})^G ; $

$$ H^0(S, \Omega^1_S) = H^0(C_1, \Omega^1_{C_1})^G \oplus H^0(C_2,
\Omega^1_{C_2})^G;
$$
$$ H^1(C_i, \mathbb{C}) = \bigoplus_{\chi \in G^*} H^1(C_i,
\mathbb{C})^{\chi} =\bigoplus_{\chi \in
G^*} (H^0(C_i, \Omega^1_{C_i})^{\chi} \oplus \overline{H^0(C_i,
\Omega^1_{C_i})^{\chi^*}}) =
$$
$$ =  \bigoplus_{\chi \in G^*} (H^0(C_i, \Omega^1_{C_i})^{\chi} \oplus H^1(C_i,
\mathcal{O}_{C_i})^{\chi});
$$
$$ H^0(S, \Omega^2_S) = \bigoplus_{\chi \in G^*} (H^0(C_1,
\Omega^1_{C_1})^{\chi} \otimes H^0(C_2,
\Omega^1_{C_2})^{\chi^*}).
$$
\end{rem}

We will use in the sequel the following quite elementary  but very 
useful result.

\begin{lemma} Let $0 \neq U' \subset U$ and $0 \neq V' \subset V$ be
complex vector spaces and set $L
:= U' \otimes V' \subset U \otimes V$. Then the subspace $L$  of 
$U\otimes V$ determines
$U'$ and $V'$.
\end{lemma}

{\em Proof.}
Let $\dashv$ denote the contraction operator $\dashv \ : (U \otimes V) 
\times V^{\vee} \ra
U.$ Then $L \dashv V^{\vee} = U'$. Similarly we get $L \dashv U^{\vee} = V'$.

\qed \\

We will see that global double Torelli holds for a huge class of
surfaces isogenous to a product of
curves, but nevertheless there are also lots of potential counterexamples. Our
first result is the following:

\begin{theo} Let $G$ be a finite abelian group acting on two curves
$C_1, C_2$  of respective genera
$g_1, g_2 \geq 2$ and  acting freely  by the product action on
$Z := C_1 \times C_2$.
Double global Torelli holds for
$S := C_1 \times C_2 /G.$, i.e., the Hodge structure of $S$ determines
the Hodge structures of $C_1$ and
$C_2$, under the following hypothesis:

A) $C_i/G $, for $i = 1,2$, has either genus $\geq2$ or it has genus 
$1$ but there is no
nontrivial subgroup $H$ of $G$ such that $C/H$ has genus $1$.
\end{theo}

{\em Proof.} Since every automorphism leaves the fundamental class of
$C_i$ invariant, we have
$H^2(C_i, \mathbb{C}) = H^2(C_i, \mathbb{C})^G = H^1(C_i,
\Omega^1_{C_i})^G =  H^1(C_i,
\Omega^1_{C_i})$. \\ 
Therefore the proof reduces to the following problem: \\
We know that $H^1(C_i, \mathbb{C}) = \bigoplus_{\chi \in G^*}
H^1(C_i, \mathbb{C})^{\chi}$ and we have
to recover the decomposition
$$ H^1(C_i, \mathbb{C})^{\chi} = H^0(C_i, \Omega^1_{C_i})^{\chi}
\oplus \overline{H^0(C_i,
\Omega^1_{C_i})}^{\chi^*}
$$
for each $\chi \in G^*$.\\ Obviously, if for $\chi \in G^*$ one of
the two above summands is zero, we
are done for this $\chi$, once we know which of the two summands is 
equal to zero.

Therefore, let
$\chi$ be a character such  that $H^0(C_i,
\Omega^1_{C_i})^{\chi} \neq 0$ and $H^0(C_i, \Omega^1_{C_i})^{\chi^*}
\neq 0$.\\
Using
$$ H^0(S, \Omega^2_S) = \bigoplus _{\chi} (H^0(C_1,
\Omega^1_{C_1})^{\chi} \otimes H^0(C_2,
\Omega^1_{C_2})^{\chi^*})
$$ and
$$ H^1(C_i, \mathbb{C})^{\chi} = H^0(C_i, \Omega^1_{C_i})^{\chi}
\oplus \overline{H^0(C_i,
\Omega^1_{C_i})}^{\chi^*}
$$
we see that the contraction
$$ H^0(S, \Omega^2_S) \dashv H^1(C_2, \mathbb{C})^{\chi}
\longrightarrow H^0(C_1,
\Omega^1_{C_1})^{\chi}
$$
is surjective, if $H^0(C_2, \Omega^1_{C_2})^{\chi{*}} \neq 0$. Using
that $g(C/G) \geq 1$ implies
$H^0(C, \Omega^1_C)^{\chi} \neq 0$ (cf. following proposition) we are
done.
\qed

\begin{prop} Let $C$ be a smooth algebraic curve of genus at least $2$
and suppose we have an effective action of  $G := \mathbb{Z} / d$
  on $C$. \\
1) Let $\chi \in G^*$: if $g(C/G) \geq 2$, or $g(C/G) = 1$ and $\chi$ is a primitive 
character $(mod \ d)$,
then $H^0(C, \Omega^1_C)^{\chi} \neq 0$.\\
2)  Let $G_{\chi}$ be the subgroup which is the kernel of
a character $\chi$ and suppose that $g(C/ G_{\chi}) \geq 2 $: if $g(C/G) = 1$, then $H^0(C, \Omega^1_C)^{\chi} \neq 0 $. If
$g(C/G) = 0$ then
either $H^0(C, \Omega^1_C)^{\chi} \neq 0 $ or $H^0(C, 
\Omega^1_C)^{\chi*} \neq 0.$

Clearly, if $\chi \neq 0$, and $ g(C/ G_{\chi}) = 0$, then 
$H^0(C, \Omega^1_C)^{\chi}
= H^0(C, \Omega^1_C)^{\chi*} =0$. If instead $g(C/G_{\chi}) = 1$ and $|G / G_{\chi}| \geq 3$, either $H^0(C, \Omega^1_C)^{\chi} = 0 $ or $H^0(C, 
\Omega^1_C)^{\chi*} = 0.$

3) If
$g(C/G) = 0$, then, assuming that the number of branch points is
at least $4$, the cardinality of the set
$\{ \chi \in G^* : H^0(C, \Omega^1_C)^{\chi} =0\}$ is strictly smaller than
$\frac{d-1}{2}$.

\end{prop}

{\em Proof.} 1) Let $P_1, \ldots , P_r \in Y$ be the branch points of
the map $C \rightarrow C/G =: Y$.

After we fix a generator for the group $G$, thus also for the group 
of characters
$G^*$, each $P_i$ determines an isotropy subgroup $H_i \cong \frac{d}{d'_i}
\ 
\ZZ /  d\ZZ \subset
\ZZ / d \ZZ$ and a character of the representation of $H_i$ on the 
tangent space $
T_{C,Q_i}$, $Q_i$ being any point lying over $P_i$.

As we shall see later, representing this character by an integer $a_i 
\ (mod\  d'_i)$,
with $( a_i , d'_i) = 1$, we obtain $m'_i$ by the equation $ a_i m'_i 
\equiv 1 (mod \
d')$. We then obtain
  natural numbers $1 \leq m_j \leq d - 1$ setting $ m_j:  = m'_j \frac{d}{d'_j}$
and we have
$$ dL \equiv \sum m_j P_j =: \hat{B}.
$$
Let $\sigma_j $ be the unique section of $\hol_Y(P_j)$ vanishing on 
$P_j$: then
the equation $z^d = \prod
\sigma_j^{m_j}$ in
$\mathbb{L}$ defines a  singular covering $p: X \rightarrow
Y$ (such that $C$ is the normalization of $X$). We have
$$ p_* \omega_X = \omega_Y \oplus \omega_Y(L) \oplus \ldots \oplus
\omega_Y(L^{d-1}).
$$ A local generator of $\omega_Y(L^i)$ is given by $dx \cdot z^{-i}$
and at $P_j$ we have (setting $m
:= m_j$, $x = \sigma_j$) the local equation $z^d = x^m$ for $X 
\subset \mathbb{L}$.

  We will investigate now when is
$\varphi (x) dx \cdot z^{-i}$ regular on $C$. \\ Let $r$ be the
greatest common divisor of $d$ and $m$,
and write $m = rm'$, $d = rd'$. Then the equation $z^d = x^m$ decomposes in
$$
\prod _{\epsilon^r = 1} (z^{d'} - \epsilon x^{m'}) = 0.
$$

We choose a point on the normalization $C$ of $X$ and let $t$ be a
local coordinate of $C$ around
this point. Then the cyclic group $H_j \cong \{ \zeta | \zeta^{d'} = 
1\}$ acts locally
by sending
$$ z \ra \zeta z , \  t \ra \zeta^a t . $$
Moreover, the
equations
$$ z = t^{m'} \epsilon ^{\frac{1}{d'}},
$$
$$ x = t^{d'}
$$ give a parametrisation of the branch of $C$ over $P_j$. Then
$\varphi (x) dx \cdot z^{-i}$ (on
$Y$) pulls back to
$$
\varphi (t^{d'}) d(t^{d'}) (t^{m'} \epsilon^{\frac{1}{d'}})^{-i} = d'
t^{d'-1}\varphi (t^{d'})dt \cdot
t^{-im'} \cdot \epsilon ^{\frac{-i}{d'}}).
$$ This is regular iff $\varphi (x) x \cdot x^{- \frac{im' + 1}{d'}}$
has order at least $0$, i.e. iff
$ord \varphi + 1 \geq \frac{im' + 1}{d'}$, or, equivalently, $\lceil 
y \rceil : = - [ -y]$
denoting the round up of a real number $y$, iff
$ord
\varphi \geq \lceil \frac{im' +
1}{d'} \rceil - 1 = \lceil \frac{im + r}{d} \rceil - 1$. \\
Therefore we have been able to compute the eigenspace for the $i-th$ 
character as:
$$ H^0(C; \Omega^1_C)^i = H^0(Y,\mathcal{L}_i),
$$ 
where
$$
\mathcal{L}_i = \omega_Y (iL - \sum_j (\lceil \frac{im_j + r_j}{d} \rceil - 1))
$$
is a line bundle of degree
$$ 2g(Y) - 2 + \sum_j (\frac{im_j}{d} - (\lceil \frac{im_j + r_j}{d}
\rceil - 1)).
$$ We remark that
$$
\frac{im_j}{d} - (\lceil \frac{im_j + r_j}{d}
\rceil - 1) =\frac{im'_j}{d'} - (\lceil \frac{im'_j + 1}{d'} \rceil - 1) =
\frac{im'_j}{d'} + 1 -
([\frac{im'_j}{d'}] + 1) = \{\frac{im'_j}{d'}\} .
$$
I.e., we obtain that the summands provide exactly the fractionary part
$\{\frac{im_j}{d}\}$ of $\frac{im_j}{d}$, in other words, the 
remainder class of $i m_j (mod
\ d)$, divided by $d$.

Therefore if $g(Y) \geq 2$ we see that $deg \mathcal{L}_i \geq 2g(Y)
- 2 \geq 2$, whence we get $h^0(C,
\Omega^1_C) ^i \geq g(Y) -1 \neq 0$.

If instead $g(Y) = 1$ and $\chi$ is a primitive character, then the 
G.C.D.  between $i$
and $d$ is $1$, whence $\{\frac{im_j}{d}\} > 0$ for each $m_j \neq 0$:
we conclude since, the genus of $C$ being $\geq 2$, there is at least 
one $m_j \neq 0$.\\

2) If $\chi$ is a non primitive character, we simply observe that 
$\chi$ yields a
primitive character $\chi'$ of $ G / G_{\chi}$, and $ H^0(C, 
\Omega^1_C)^{\chi}$ is the
pull back of  $H^0(C/ G_{\chi}, \Omega^1_C)^{\chi'}$, so that the 
first assertion is a
direct consequence of 1).

Assume now that $g(C/ G_{\chi}) \geq 2$ and  $g(C/G) = 0$. We apply 
the basic estimate
we used in the proof of 1):  i.e.,  we have a primitive character and then
  $ H^0(C, \Omega^1_C)^{\chi} \neq 0$ if $ \sum_j \{\frac{im_j}{d}\} \geq 2.$

Without loss of generality, since we have a primitive character, we 
may assume $i=1$.
$ H^0(C, \Omega^1_C)^{\chi} \neq 0$ unless $\sum_j m_j = d$. But our 
assertion  holds
since the dual character corresponds to $d-1$, and  $ H^0(C, 
\Omega^1_C)^{\chi*} \neq 0$
unless $\sum_j (d -m_j) = d$. This is a contradiction, since then
$ 2 d = \sum_j (d -m_j + m_j) = d( \sum_j 1)$, there are exactly 
$r=2$ branch points,
whence $g(C/ G_{\chi}) = 0$.\\

3) Assume finally  $g(Y) = 0$, i.e. $Y = \mathbb{P}^1$. Obviously $H^0(C,
\Omega^1_C)^0 = H^0(Y,
\Omega^1_Y) =0$. For $i > 0$ we have, by our previous calculation:
$$ H^0(C, \Omega^1_C)^i \neq 0
$$ if and only if
$$
\sum _j \{ \frac{im_j}{d} \} \geq 2.
$$ Then the assertion follows from the following lemma.
\qed

\begin{lemma} Let $r, d, m_1, \ldots m_r$ be natural numbers such that
$r \geq 4$ and $m_1 + \ldots m_r
= d$. We define for each $1 \leq i \leq d - 1$ a natural number
$\lambda (i)$ by the equality
$$
\overline{im_1} + \ldots \overline{im_r} = \lambda (i) d,
$$
where $0 \leq \overline{z} \leq d-1$ denotes the remainder modulo $d$ of a natural
number $z$. Then
$$
\sharp \{i : \lambda (i) = 1 \} < \frac{d-1}{2}.
$$
\end{lemma}

{\em Proof.} Assume that $r > 4$. Then we set $M_{r-1} := m_{r-1} +
m_r$ and it follows
$$
\lambda ^r (i) d = \overline{im_1} + \ldots \overline{im_r} \geq
\overline{im_1} + \ldots
\overline{iM_{r-1}} = \lambda ^{r-1} (i) d.
$$ 
Therefore without restriction we can assume $r = 4$. For $d = 4$ we
have $m_1 = m_2 = m_3 = m_4 =1$
and the claim is obvious. Let now $d + 1 > 4$ and $d + 1 = m_1 +
\ldots + m'_4$. Without restriction we
can assume that $m'_4 > 1$, hence we write $m'_4 = m_4 + 1$. Then
$m_1 + \ldots + m_4 = d$ and by
induction we get $\sum _ j \overline{im_j}^d = \lambda '(i)d$ and
$\sharp \{i : \lambda' (i) = 1 \} <
\frac{d-1}{2}$.\\ 
We write for $1 \leq j \leq 3$:
$$ im_j = a_jd + \overline{im_j}^d = b_j (d+1) +\overline{im_j}^{d+1}
$$ and
$$ i(m_4 + 1) = \gamma d + \overline{im_4}^d + i = c (d+1)
+\overline{i(m_4 + 1)}^{d+1}.
$$ Then
$$
\overline{i(m_4 + 1)}^{d+1} = \overline{im_4}^d + (\gamma - c) d + i - c,
$$ whence
$$
\lambda (i) (d+1) = \overline{im_1}^{d+1} + \ldots + \overline{i (m_4 
+1)}^{d+1} =
$$
$$ = \overline{im_1}^d + \ldots + \overline{im_4}^d + (a_1 - b_1) d +
\ldots + (a_3 - b_3) d + (\gamma
- c) d + i - b_1 - b_2 - b_3 -c.
$$ We remark that
$$ i d = (a_1 + a_2 + a_3 + \gamma) d + \lambda '(i) d.
$$ Therefore $\lambda '(i) = 1$ implies that $a_1 + a_2 + a_3 +
\gamma = i - 1$. Analogously $\lambda
(i) = 1$ implies that $b_1 + b_2 + b_3 + c = i - 1$. \\

Assume that $\lambda (i) = 1$: then
$$
\lambda '(i) d + (a_1 - b_1) d + \ldots + (\gamma - c) d + 1 = d + 1.
$$ Since $a_1 - b_1, \ldots, \gamma - c \geq 0$ we see that the above
equality implies that $\lambda '(i) =
1$. Therefore $\sharp \{ i \leq d-1 : \lambda (i) = 1\} <
\frac{d-1}{2}$ and we are done if we show that $\lambda (d) \neq 1$.

Since
$$
\lambda (d) (d + 1) = \overline{dm_1}^{d+1} + \ldots + \overline{d
(m_4 +1)}^{d+1},
$$
if $\lambda (d) = 1$ we get
$$
2 (d + 1) = m_1 + \ldots + (m_4 + 1) + \overline{dm_1}^{d+1} + \ldots + 
\overline{d
(m_4 +1)}^{d+1} =
$$
$$
= (m_1 + \overline{dm_1}^{d+1}) + \ldots + ((m_4 + 1) + 
\overline{d(m_4 + 1)}^{d+1}) =
4(d+1),
$$
which is absurd.

\qed

\section{Another series of examples}

In the following we will construct a series of examples of surfaces
isogenous to a higher product, for which we are not able to prove 
global double Torelli
using theorem 4.4 or proposition 4.5.
It turns out nevertheless that for these surfaces the local period 
map is injective. \\

It seems therefore an interesting problem to decide whether global double
Torelli holds for this class of surfaces (or for similar classes).

We consider for
any natural number $d \geq 3$ the
group $G = G_d = \mathbb{Z} / d \oplus \mathbb{Z} / d$. Let $\pi :
C_1 \longrightarrow C_1 / G =
\mathbb{P}^1$ be the covering branched in $P$, $P'$ (with local 
monodromy given by
$(1,0)$, respectively $(-1,0)$) and in $P_1,
\ldots P_d$ (with  $(0,1)$ as local monodromy).

  $C_1$ is then the curve in the weighted projective plane 
$\mathbb{P}(1,1,d)$
with coordinates $(x_0, x_1, z)$ defined by the equation (homogeneous 
of degree $d^2$)
$$ z^d = \prod_{j=1}^d (x_1^d - \alpha _j x_0^d).
$$
$\pi$ is given by $\pi (x_0, x_1, z) : = (x_0^d, x_1^d)$. We set for 
convenience
$x: = x_1/ x_0$, $u := x^d$, $u_0 := x_0^d, u_1 :=  x_1^d.$

Then the canonical sheaf of $C_1$ is $\hol_{C_1} (d^2 -d-2)$ and it 
is easy to see that,
in affine coordinates, the space of holomorphic $1$-forms can be 
written as follows,
where the $F_j(x)$'s are polynomials of degree $\leq j$:
$$ H^0(C_1, \Omega^1_{C_1}) = \{ F_{d-2}(x) \frac{dx}{z} + F_{2d-2}(x)
\frac{dx}{z^2} + \ldots + F_{d(d-1)-2}(x)
\frac{dx}{z^{d-1}} \}.
$$
Accordingly, the canonical map is given by $\Phi(x_0, x_1, z) = (z^i 
x_0^h x_1^r)$,
where $ 0 \leq i \leq d-2$, $ h + r + d i = (d^2 -d-2)$.

We can now write the cohomology table of $C_1$, where we write in the
place $(a,b)$ the dimension of
the eigenspace of $H^0(C_1, \Omega^1_{C_1})$ belonging to the
character $(a,b)$, i.e. $dim H^0(C_1,
\Omega^1_{C_1})^{(a,b)}$, noting that $\frac{dx \cdot x^i}{z^j}$
belongs to the character $(-j,
\overline{i+1})$.

\medskip

\begin{tabular}{|l|c|c|c|c|c|c|}
\hline &$a = 0$ &
$a = 1$ &
$a = 2$ &
$a = 3$ &
$ \ldots $ &
$a = d-1$
\\
\hline \hline
$b = 0$ & $0$ & $0$ & $1$ & $2$ & $ \ldots $ & $d-2$\\
\hline
$b = 1$ & $0$ & $1$ & $2$ & $3$ & $ \ldots $ & $d-1$\\
\hline
$\ldots$ &  $\ldots $ & $ \ldots $ & $\ldots $ & $ \ldots $ & $
\ldots $ & $\ldots $ \\
\hline
$b = d-2$ & $0$ & $1$ & $2$ & $3$ & $ \ldots $ & $d-1$\\
\hline
$b = d-1$ & $0$ & $1$ & $2$ & $3$ & $ \ldots $ & $d-1$\\
\hline
\end{tabular}

\medskip We take now the covering $\varphi : C_2 \longrightarrow C_2
/ G = \mathbb{P}^1$ branched in
$P$, $P'$ (with local monodromies $(1,0)$ and $(-1,0)$) and in $Q, Q'$
(with local monodromies $(0,1)$, $(0,-1)$).

We may assume without loss of generality that the points $P$, $P'$
are the respective points $ v =0, v = \infty$, whereas the points
$Q, Q'$ are the respective points $ v =1, v = - \lambda$.

We see easily that $C_2 \subset \mathbb{P}^1 \times \mathbb{P}^1$ is given by
the equation
$$ w_1^d (y_1^d + \lambda y_0^d) = w_0^d (y_1^d - y_0^d),
$$
$\varphi ((y_0,y_1)(w_0, w_1)) : = (y_0^d,y_1^d)$ and in affine coordinates
$C_2$ is the fibre product of two cyclic coverings:
$$ y^d = v,
$$
$$ w^d = \frac{v-1}{v+1}.
$$
We see  immediately that
$$ H^0(C_2, \Omega^1_{C_2}) = \{ F_{\leq d-2, \leq d-2} (w,y)
\frac{dy}{w^{d-1} (y^d + 1)} \} =  \{
F_{\leq d-2, \leq d-2} (w,y) \frac{wdy}{v-1} \}.
$$
Therefore we have the following cohomology table:

\medskip

\begin{tabular}{|l|c|c|c|c|c|c|}
\hline &$a = 0$ &
$a = 1$ &
$a = 2$ &
$a = 3$ &
$ \ldots $ &
$a = d-1$
\\
\hline \hline
$b = 0$ & $0$ & $0$ & $0$ & $0$ & $ \ldots $ & $0$\\
\hline
$b = 1$ & $0$ & $1$ & $1$ & $1$ & $ \ldots $ & $1$\\
\hline
$\ldots$ &  $\ldots $ & $ \ldots $ & $\ldots $ & $ \ldots $ & $
\ldots $ & $\ldots $ \\
\hline
$b = d-2$ & $0$ & $1$ & $1$ & $1$ & $ \ldots $ & $1$\\
\hline
$b = d-1$ & $0$ & $1$ & $1$ & $1$ & $ \ldots $ & $1$\\
\hline
\end{tabular}

\bigskip In order to have $G$ operating freely on the product $C_1
\times C_2$ we have to twist the
action of $G$ on $C_2$. We assume for simplicity that $d$ is a prime number.
Let now $r \neq 0$ and $\neq 1$. Then twisting the action of $C_2$ by 
the automorphism of
$\mathbb{Z} / d \mathbb{Z} \oplus \mathbb{Z} / d \mathbb{Z}$ given by
$$
(1, 0) \mapsto (1, 1), \ \ \ (0, 1) \mapsto (r, 1),
$$
we see that the stabilizers of the twisted action on $C_2$ are now $<(1,
1)>$  and $<(r, 1)>$, whence $G$ acts now freely on $C_1 \times C_2$ and 
the cohomology table of $C_2$ becomes now:

\medskip

\begin{tabular}{|l|c|c|c|c|c|c|}
\hline &$a = 0$ &
$a = 1$ &
$a = 2$ &
$a = 3$ &
$ \ldots $ &
$a = d-1$
\\
\hline \hline
$b = 0$ & $0$ & $1$ & $1$ & $1$ & $ \ldots $ & $1$\\
\hline
$b = 1$ & $1$ & $0$ & $*$ & $*$ & $ \ldots $ & $*$\\
\hline
$\ldots$ &  $\ldots $ & $ \ldots $ & $\ldots $ & $ \ldots $ & $
\ldots $ & $\ldots $ \\
\hline
$b = d-2$ & $1$ & $*$ & $*$ & $*$ & $ \ldots $ & $*$\\
\hline
$b = d-1$ & $1$ & $*$ & $*$ & $*$ & $ \ldots $ & $0$\\
\hline
\end{tabular}

\bigskip
I.e., the diagonal is zero and the remaining zeroes are at the places 
$(nr, n)$,
where $0 \leq n \leq d-1$. In particular we see that e.g.
$$
H^0(C_1, \Omega_{C_1}^1)^{(1,1)} \neq 0 \neq H^0(C_1, 
\Omega_{C_1}^1)^{(d-1, d-1)}
$$
whereas
$$
H^0(C_2, \Omega_{C_2}^1)^{(1,1)} = 0 = H^0(C_2, \Omega_{C_2}^1)^{(d-1, d-1)}.
$$
Therefore we cannot reconstruct $H^0(C_1, \Omega_{C_1}^1)^{(1,1)}$ as well as
  $H^0(C_1, \Omega_{C_1}^1)^{(d-1, d-1)}$ from the Hodge structure of $S$.\\

By theorem \ref{thm1} we can now easily verify that the infinitesimal 
Torelli map is injective.
In fact, $H^0(C_i, \hol_{C_i}( 2 K_{C_i}))^{G} =H^0(\mathbb{P}^1, 
\hol_{\mathbb{P}^1}(-4 + r_i)$, where
$r_i$ is the number of branch points of $C_i \ra \mathbb{P}^1$.

For $i=2$ we get a space of dimension $1$, therefore it suffices to 
observe that there is
a non zero summand in
$$
\bigoplus_{\chi \in G^* : H^0(C_1,\Omega^1_{C_1})^{\chi} \ \neq 0} H^0(C_2,
\Omega^1_{C_2})^{\chi} \otimes H^0(C_2, \Omega^1_{C_2})^{\chi^*}.
$$
For $i=1$ we must obtain all monomials of degree $d-2$ in $(u_0, 
u_1)$, and for this purpose
it suffices, since the pairing is non degenerate, to find a character space
$H^0(C_1,\Omega^1_{C_1})^{\chi}$ of dimension $d-1$ such that the 
character space
$H^0(C_1,\Omega^1_{C_1})^{\chi*}$ is non zero, and such that likewise 
$H^0(C_2,\Omega^1_{C_2})^{\chi}
\neq 0$. \\

We omit here to prove the following

\begin{prop}
The canonical system of $S$ has base points, and the canonical map is 
birational onto its image.
\end{prop}

\begin{footnotesize}
\noindent
{\bf Acknowledgement .}
We would like to thank Gerard van der Geer for providing numerical evidence
for Lemma 4.6 when we conjectured it.
\end{footnotesize}

\end{document}